\documentclass[journal]{IEEEtran}
\usepackage[utf8]{inputenc}
\usepackage{textgreek}

\IEEEoverridecommandlockouts                              

\usepackage{graphics} 
\usepackage{epsfig} 
\usepackage{mathptmx} 
\usepackage{times} 
\usepackage{amsmath} 
\usepackage{amssymb}  
\usepackage{cleveref}
\usepackage{algorithm}
\usepackage{algpseudocode}
\usepackage{graphicx}
\usepackage{subcaption}
\usepackage{xcolor}
\usepackage{tikz}
\usetikzlibrary{arrows.meta, positioning, shapes.geometric}
\usepackage{csquotes}
\usepackage{cite}

\title{\LARGE \bf
Collaborative Decision-Making and Optimal Utilization of Pathfinding Flights during Convective Weather
}

\author{Jimin Choi$^{1}$, Husni R. Idris$^{2}$, {Huy T. Tran}$^{3}$, {Max Z. Li}$^{4}$%
\thanks{*This work was supported by NASA's Transformative Tools and Technologies Project, Award No. 80NSSC23M0221.}%
\thanks{$^{1}$Jimin Choi is with the Department of Aerospace Engineering, University of Michigan,
        Ann Arbor, MI 48109, USA, 
        {\tt\small jiminch@umich.edu}}%
\thanks{$^{2}$Husni R. Idris is with the Aviation Systems Division, NASA Ames Research Center, Moffett Field, CA 94035, USA, {\tt\small husni.r.idris@nasa.gov}}%
\thanks{$^{3}$Huy T. Tran is with the Grainger College of Engineering,
Department of Aerospace Engineering, University of Illinois Urbana-Champaign, Champaign, IL 61820, USA,
        {\tt\small huytran1@illinois.edu}}%
\thanks{$^{4}$Max Z. Li is with the Departments of Aerospace Engineering, Civil and Environmental Engineering, Industrial and Operations Engineering, University of Michigan,
        Ann Arbor, MI 48109, USA, 
        {\tt\small maxzli@umich.edu}}%
}

\begin{document}

\maketitle
\thispagestyle{empty}
\pagestyle{empty}


\begin{abstract}
Air traffic operations are strongly influenced by convective weather, and one common response is \emph{pathfinder operations}, in which a designated aircraft tests the viability of weather-impacted airspace and routes. Despite relatively routine use in practice, how pathfinder operations evolve under uncertainty and how the pathfinder decision-making process unfolds are largely treated as exogenous. Addressing this gap requires jointly modeling weather-driven system accessibility, flight responses to pathfinder offers, and the sequencing of those offers to improve outcomes. We develop a unified analytical framework that connects weather-driven system state transitions, flight acceptance decisions, and the sequencing of pathfinder offers. We first construct a four-state Markov chain to model stochastic closure and reopening of exit points, or \emph{fixes}, out of the terminal departure airspace surrounding a major airport, pathfinder selection, and pathfinding execution, and analyze its steady-state behavior to characterize long-term capacity and delay implications. We introduce utility-based decision models for flights, air traffic control (ATC), and dispatchers, and analyze worst-case collective rejection to quantify system vulnerability under selfless behavior and uncertainty. Finally, we formulate optimization problems that model ATC-initiated and dispatcher-initated pathfinder offers, with the goal of optimizing the \emph{sequence} of pathfinder offers. Using a discrete event simulation for a major US airport, we show that ATC- and dispatcher-driven objectives lead to distinct, near real-time sequencing strategies, providing the first formal decision models for pathfinder operations under weather uncertainty.
\end{abstract}


\section{Introduction}\label{sec:introduction}
\subsection{Background and Motivation}
\label{sec:background_motivation}
Convective weather is one of the dominant sources of disruption in commercial aviation, affecting both flight operations and air traffic management (ATM). To mitigate safety risks and avoid exposure to severe turbulence, thunderstorms, and hazardous atmospheric conditions, airlines routinely reroute flights away from regions of strong convective activity \cite{gultepe2019review}. In response to weather hazards, air traffic control (ATC) imposes temporary closures that impact specific terminal and en route airspace \cite{mitchell2006airspace}. When such restrictions persist, they can significantly degrade airline schedule reliability and network performance \cite{aerospace10030288}. The resulting economic burden is substantial, as illustrated by total delay-related costs in the United States exceeding \$33 billion in 2019 alone \cite{FAA2022}.

Weather-driven airspace restrictions have particularly severe operational consequences in high-density terminal environments such as the New York terminal area (N90 TRACON) \cite{allan2001nydelay}. Aircraft access to and from the terminal region is constrained to a set of departure fixes defined by fixed geographic coordinates and altitudes. When these fix locations, along with jet routes emanating from these locations, are impacted by convective weather, they are often closed by ATC to reduce deviations and excess radar vectoring \cite{interview}. When multiple departure fixes are unavailable, airport surface congestion escalates rapidly, leading to extended departure queues. In extreme cases, these disruptions propagate to arriving traffic through congested airport surface taxiways, lack of airport gate availability (due to departures being held at the gate), and unavailability of routes at major N90 airports \cite{mukherjee2011flight}. Moreover, empirical network-level studies have shown that such local terminal disruptions can propagate across the broader domestic and international air transportation system \cite{li2017beijing, li2021gsp}.

Timely identification of when previously restricted fixes can be safely reopened is therefore a critical operational task \cite{interview}. However, direct assessment of atmospheric conditions near these high-altitude, geographically dispersed airspace elements remains challenging, due to the inherent uncertainty and limited tactical fidelity of convective weather observation and short-term forecasting tools \cite{wolfson2006advanced, doi:10.2514/6.2011-6511}. As a result, reopening decisions often rely on exploratory flights to validate actual route usability. These flights, commonly referred to as \textit{pathfinder flights} \cite{FAA2024VP}, play a key role in reopening previously closed departure fixes, with a single successful pathfinder often yielding substantial delay reductions for downstream flights. 

\begin{figure}[htbp]
\centering
\includegraphics[width=0.65\linewidth]{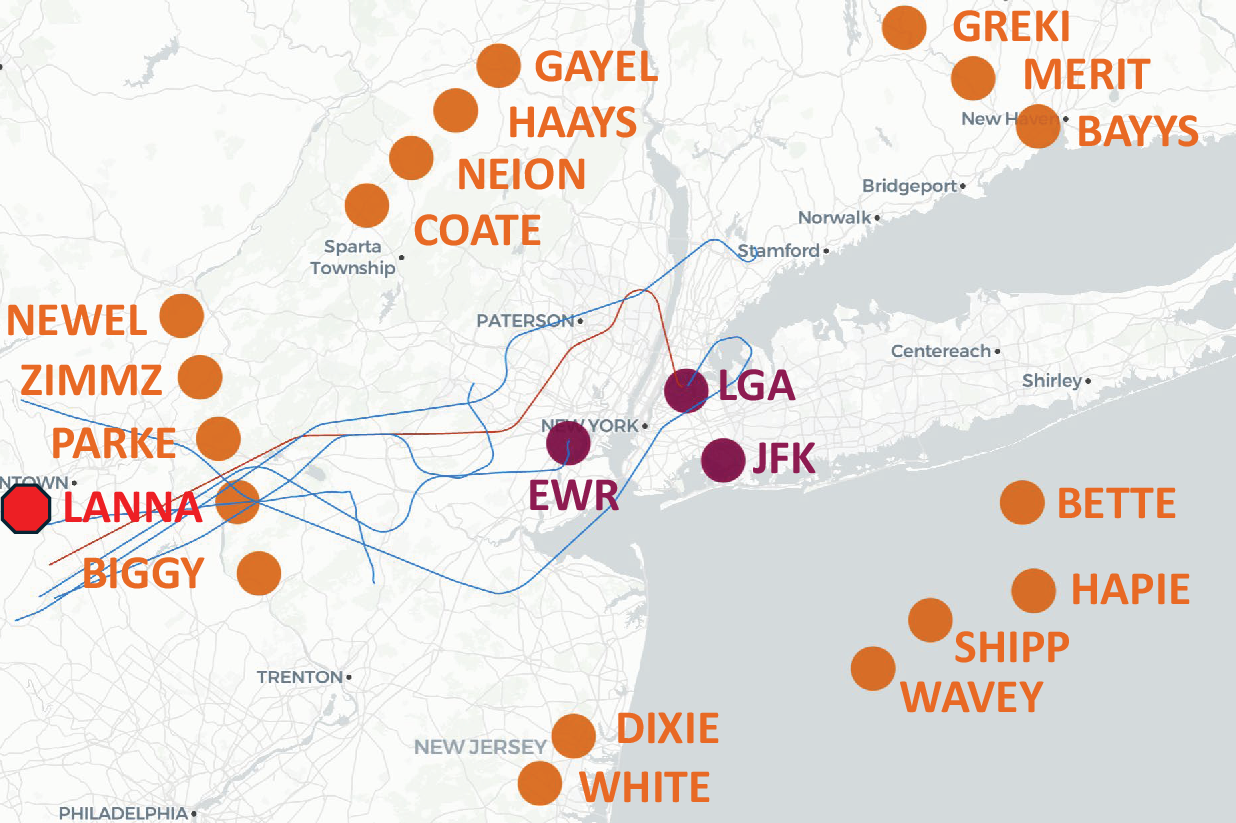}
\caption{Example of a pathfinder operation in the New York terminal area (N90). Orange circles denote departure fixes, purple circles denote major airports, and the red trajectory corresponds to the designated pathfinder flight (DAL569).}
\label{fig:pathfinder_example}
\end{figure}

\Cref{fig:pathfinder_example} shows a real-world example of a pathfinder operation in the New York terminal area identified via National Traffic Management Logs (NTMLs) from the US Federal Aviation Administration (FAA). On July~31,~2024, the designated pathfinder flight (Delta Air Lines 569, or DAL569) was requested and initiated attempted to reopen the LANNA departure fix but was rerouted in air toward PARKE due to convective weather, resulting in an unsuccessful pathfinding. 

In operational practice, when a departure fix is blocked by convective weather, ATC may request a single aircraft that is either scheduled to use that fix or willing to accept a reroute to attempt passage through the region. If the flight successfully traverses the route segment, ATC may immediately reopen the fix and clear subsequent departures through the same access point, often resulting in rapid reductions in surface congestion and delay. At the same time, the designated pathfinder aircraft faces additional operational risk and potential passenger discomfort. Despite their operational importance, pathfinder selection and execution remain largely unmodeled. This motivates analytical frameworks that model pathfinder decision-making and associated operational trade-offs.

 \subsection{Related Work}
\label{sec:related_work}
\subsubsection{Weather Impacts and Terminal Airspace Congestion}
Extensive prior work has examined how convective weather degrades airspace capacity and disrupts traffic flows, motivating the development of weather-aware traffic flow management (TFM) strategies. Foundational studies showed how thunderstorms reduce sector capacity and drive network-wide congestion, leading to increased delays and cancellations \cite{mitchell2006airspace,gultepe2019review, aerospace10030288}. Subsequent research has developed probabilistic and optimization-based approaches for modeling weather-impacted routing, capacity uncertainty, and delay propagation across the airspace system \cite{pfeil2012identification, LUI2022103811,yang2018strategic4d,nilim2004stochastic}. Building on this system perspective, a growing body of work has focused specifically on terminal airspace operations, examining how convective weather leads to nonlinear amplification of surface congestion, departure delays, and arrival holding, particularly at high-density hubs such as the New York airports \cite{mukherjee2011flight,li2017beijing}. These results underscore that managing access to departure and arrival fixes under convective weather is central to system performance. Related work has also explored flight procedural design \cite{fowler2012getting}, integrated surface and airspace operations \cite{badrinath2019integrated}, and the design of alternative airspace configurations and geometries \cite{granberg2019framework}. Despite establishing that convective weather and fix availability are central drivers of terminal congestion, this body of work treats fix usability as exogenous and does not model the decision-driven probing and reopening of fixes through pathfinder operations.

\subsubsection{Operational Traffic Flow Management Tool}
On the operational side, the FAA has deployed a number of TFM programs to mitigate weather-driven disruptions. Collaborative initiatives such as the Collaborative Trajectory Options Program (CTOP) enable airlines to submit routing options when demand is expected to exceed capacity and to manage flows around convective weather more flexibly \cite{faa2014ac90115, kulkarni2018ctop}. The Route Availability Planning Tool (RAPT) supports strategic and tactical planning by forecasting the temporal availability of departure fixes based on convective forecasts \cite{robinson2009rapt}. These programs illustrate that deciding when and how to route traffic through weather-impacted structures is recognized as an important operational problem. However, these tools treat route or fix availability as an external input inferred from weather products, rather than as the outcome of active assessment by designated pathfinder flights. As a result, existing operational TFM tools do not model who should serve as a pathfinder, how acceptance decisions are formed, or how the sequence of pathfinder offers should be constructed under uncertainty.

\subsubsection{Optimization and Control in Air Traffic Management}
In parallel, a rich literature has developed decision-making and optimization tools for air traffic management. Stochastic optimization and Markov decision process formulations have been applied to departure scheduling, surface traffic management, and flow control under uncertainty, often focusing on minimizing expected delay or maximizing throughput \cite{mukherjee2007ground, mukherjee2009reroute, scala2020dss}. Queueing and control-theoretic models of the departure process have been used to design and evaluate dynamic pushback and departure metering strategies \cite{simaiakis2012dynamic, simaiakis2016queuing}. Beyond centralized optimization, game-theoretic and agent-based models have been used to study strategic interactions between airlines and ATC and to represent heterogeneous decision-making at the flight level \cite{evans2016airline, bongiorno2017abm}. More recently, learning-based approaches such as deep and offline reinforcement learning have been applied to airport departure metering and runway configuration decisions \cite{ali2022drl, memarzadeh2025runway}, while data-driven trajectory prediction models integrate operational and weather information to forecast 4D flight trajectories \cite{shi2021lstm4d}. While these approaches provide powerful tools for centralized and decentralized traffic management under uncertainty, they generally assume that airspace availability is given, and do not capture exploratory, sequential, and incentive-driven decision processes such as pathfinder selection and acceptance.

\subsubsection{Pathfinder Operations in Prior Literature}
Pathfinder-like concepts have primarily appeared in operational documentation. FAA NTML logs and briefing materials describe the use of pathfinder flights to test whether weather-affected departure fixes can be used safely, and recent FAA working groups have focused on pathfinder procedures and coordination mechanisms \cite{FAA2024Pathfinder}. One academic study discussed the role of exploratory flights in evaluating route usability under convective weather and their influence on traffic flow decisions \cite{Weber2007}. However, these works remain largely descriptive or procedural and do not provide a formal decision-theoretic model of pathfinder operations.

\subsection{Research Gap and Our Contributions}
\label{sec:contributions}


Convective weather is a dominant source of disruption in commercial aviation, and a wide range of optimization, queueing, and learning-based tools have been developed to manage weather-impacted traffic. In practice, pathfinder operations are also routinely used to assess airspace and route usability in real time. However, existing models treat fix usability as an exogenous weather-driven state and lack a unified framework that captures pathfinder operations as an endogenous stochastic, multi-agent decision process. In particular, prior work does not model voluntary, utility-driven acceptance decisions by individual flights, the resulting risk of collective rejection, or the sequential allocation of pathfinder offers by ATC and airline dispatchers. Among these, the pathfinder selection stage is uniquely governed by stakeholder decisions rather than weather dynamics, making it the primary control lever for decision support and optimization. This motivates a unified modeling framework that links stochastic fix availability, flight acceptance behavior, and sequential offer allocation for quantitative analysis under operational uncertainty.

This paper establishes the first analytical framework that formalizes pathfinder operations as a coupled stochastic and multi-agent decision-making problem under convective weather. Our contributions are threefold. First, we develop a Markov chain model that characterizes the long-run availability of weather-impacted departure access points and their implications for capacity and delay. Second, we introduce utility-based models of flight acceptance and ATC and airline dispatcher decision-making. We then derive worst-case collective rejection bounds under selfless behavior and environmental uncertainty. Finally, we formulate and solve sequencing optimization problems under ATC- and dispatcher-initiated control, and show that their objectives yield structurally different offer sequences. A preliminary version of this work was presented at the 2025 IEEE International Conference on Intelligent Transportation Systems (ITSC) \cite{choi2025pathfinders}. This manuscript substantially extends the previous version by incorporating stakeholder interviews to refine the utility models, developing a Markov chain–based operational analysis, and introducing optimization formulations and numerical results for sequencing pathfinder offers.

\section{Methodology}\label{sec:methodology}
This section presents a unified modeling and decision framework for pathfinder operations under convective weather. We first introduce a Markov chain model that captures the stochastic evolution of fix closure, selection, execution, and reopening, providing a system characterization of long-run capacity and delay. We then develop utility-based decision models for flights, ATC, and airline dispatchers to formalize acceptance behavior and stakeholder objectives, and analyze worst-case collective rejection to quantify system vulnerability under self-interested, selfless, and uncertain decision environments. Finally, we formulate optimization problems that determine the sequence of pathfinder offers under ATC- and dispatcher-initiated control, linking individual acceptance behavior to system performance.

\subsection{System Operational Dynamics via Markov Chains}
\label{sec:modeling_markov_chain}
The temporal dynamics of pathfinder operations are represented using a Markov chain model that captures uncertain transitions between key operational states for a specific departure fix. Steady-state analysis is used to characterize the resulting long-term system behavior.
\begin{figure}[htbp]
  \centering
  \includegraphics[width=\linewidth]{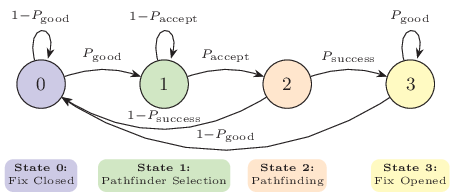}
  \caption{Markov chain representation of the aircraft \mbox{pathfinding} process. Transition probabilities are governed by \(P_{\mathrm{good}}\), \(P_{\mathrm{accept}}\), and \(P_{\mathrm{success}}\), capturing weather observations, acceptance of pathfinder offers, and pathfinding success, respectively.}
  \label{fig:markov_chain}
\end{figure}
A Markov chain describes a stochastic system in which state transitions are governed exclusively by the current state, without explicit dependence on the sequence of past states. This memoryless structure makes the framework particularly effective for representing decision-making processes under uncertainty and is widely used in areas such as reliability and systems analysis \cite{markov_chain}. Building on this modeling paradigm, we represent the evolution of the pathfinding process under convective weather using a discrete-time Markov chain, in which probabilistic transitions occur between a set of operational states. The resulting state-transition structure is shown in \Cref{fig:markov_chain}.

The proposed Markov chain is parameterized by three transition probabilities $P_\mathrm{good}$, $P_\mathrm{accept}$, $P_\mathrm{success} \in [0,1]$, each representing the probability of a well-defined operational event. $P_\mathrm{good}$ characterizes the likelihood that weather conditions are sufficiently favorable for the system to initiate pathfinding rather than remain in a closed state. $P_\mathrm{accept}$ denotes the conditional probability that a candidate aircraft agrees to undertake the pathfinder attempt when requested. Finally, $P_\mathrm{success}$ represents the conditional probability that a
pathfinding attempt results in the reopening of the affected fix.

The system evolves over four distinct operational states. State~0 (\textit{Gate Closed}) corresponds to the condition in which the departure fix is closed due to adverse convective weather. When weather observations indicate improvement with probability $P_{\mathrm{good}}$, the system transitions to State~1 (\textit{Pathfinder Selection}). In this state, a candidate flight is considered for the pathfinder. If the flight accepts the assignment with probability $P_{\mathrm{accept}}$, the process advances to State~2 (\textit{Pathfinding}). Otherwise, the system remains in State~1 and continues evaluating additional candidates. Following a pathfinding attempt, a successful outcome occurring with probability $P_{\mathrm{success}}$ moves the system to State~3 (\textit{Gate Opened}), where the fix is reopened for normal departures. An unsuccessful pathfinding forces a return to State~0. Once in State~3, the fix remains open as long as favorable weather conditions persist. If weather conditions worsen, the system returns to State~0.

The long-term dynamics of the Markov chain can be examined through steady-state analysis, which yields the limiting distribution over the four operational states introduced earlier, representing the long-run proportion of time the system spends in each state. For  simplicity, the values of $P_{\mathrm{good}}$, $P_{\mathrm{accept}}$, and $P_{\mathrm{success}}$ are treated as fixed, corresponding to a homogeneous Markov chain. Under these assumptions, the system evolution is fully characterized by the transition matrix $P$.
\begin{equation}
P = 
\begin{pmatrix}
1 - P_{\mathrm{good}} & P_{\mathrm{good}} & 0 & 0 \\
0 & 1 - P_{\mathrm{accept}} & P_{\mathrm{accept}} & 0 \\
1 - P_{\mathrm{success}} & 0 & 0 & P_{\mathrm{success}} \\
1 - P_{\mathrm{good}} & 0 & 0 & P_{\mathrm{good}}
\end{pmatrix}.
\end{equation}
Additionally, the steady-state distribution $\pi = (\pi_0, \pi_1, \pi_2, \pi_3)$ is defined as the unique probability vector satisfying the balance and normalization conditions $\pi P = \pi$ and $\sum_{i=0}^{3} \pi_i = 1$, respectively, with all transition probabilities lying in the unit interval \cite{MC_solve}. Solving this system yields the stationary occupancy of each operational state (i.e., the long-run proportion of time the system spends in each state). 

Each state represents a distinct operational condition of the affected departure fix. Departures can be released through the fix only in State~3, which is reached after a successful pathfinding attempt confirming that the fix is usable under the prevailing weather. Consequently, the stationary probability $\pi_3$ directly measures the long-run fraction of time during which departure capacity is available.

This interpretation allows the steady-state distribution to be used as a compact summary of system performance. In particular, the probability $\pi_3$ can be mapped to an effective service rate and combined with standard queueing approximations to assess throughput, stability (i.e., whether queues remain bounded over time), and delay. These operational implications of the Markov chain are examined in \Cref{sec:markov_result}.

\subsection{Stylized Modeling of Pathfinder Selection}
\label{sec:pathfinder_selection}
Building on the Markov chain model that captures system operational dynamics under convective weather, we now focus on the \textit{Pathfinder Selection} phase. Unlike the other states in the Markov chain driven by uncontrollable external factors like weather, State~1 is determined purely by stakeholder decisions and thus offers opportunity for operational improvement. This motivates a closer examination of stakeholder interactions and decision-making in this stage.

\begin{figure}[htbp]
  \centering
  \includegraphics[width=\linewidth]{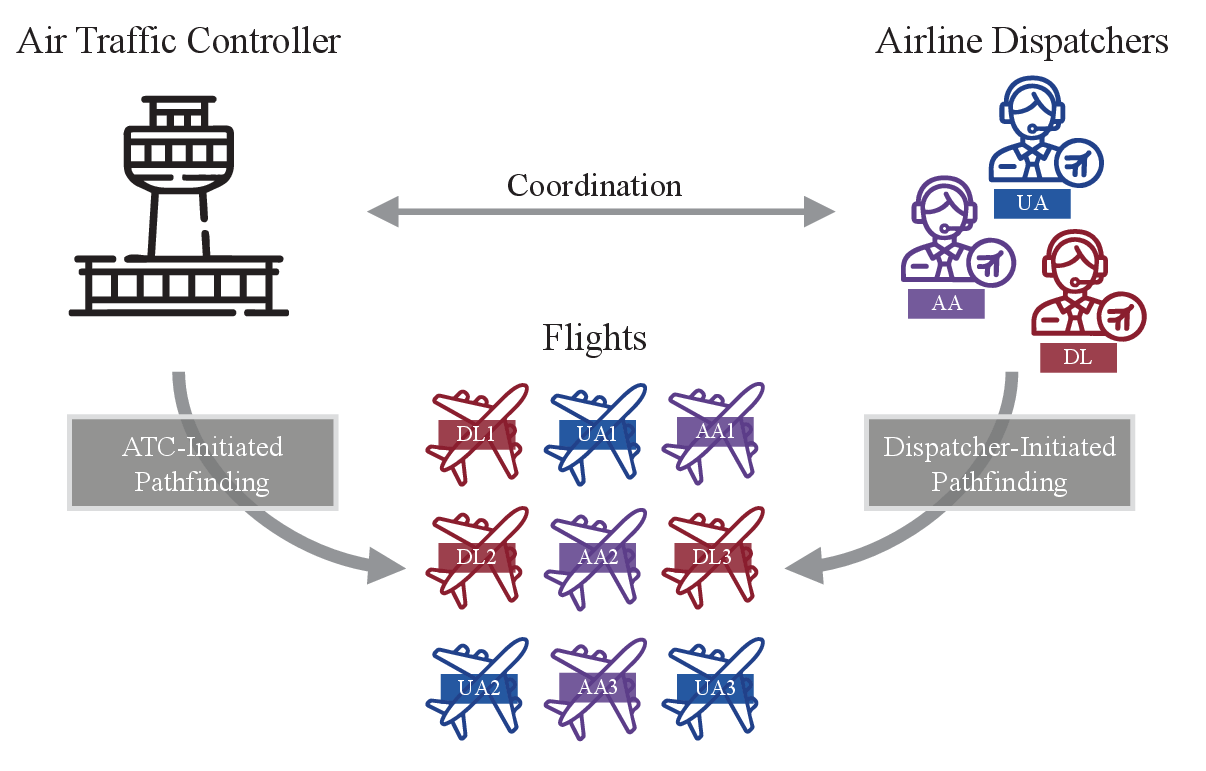}
  \caption{Stakeholders involved in the pathfinder selection. The figure shows the three primary stakeholders: ATC, airline dispatchers, and individual flights.}
  \label{fig:stakeholders}
\end{figure}

To ensure that our modeling reflects real operations, we conducted interviews with individuals involved in or affected by pathfinder activities (e.g., FAA traffic managers, airline dispatchers, and pilots) \cite{interview}. These interviews revealed the decision dynamics and information flows among the three stakeholders, summarized in \Cref{fig:stakeholders}, which inform our modeling without aiming to capture all operational contexts. \emph{ATC} manages overall traffic flow and aims to reopen closed departure fixes to restore throughput. Airline \emph{dispatchers} have similar goals, but their priorities are shaped by airline-specific factors like regulatory requirements (e.g., DOT3 rules\footnote{The U.S. Department of Transportation (DOT) three-hour tarmac delay rule limits domestic flights to a maximum of three hours on the tarmac without allowing passengers to deplane, except under safety, security, or ATC constraints \cite{dot3_rule}.}), delay costs, and downstream scheduling effects. Individual \emph{flights} decide whether to accept a pathfinder offer, and their decisions determine whether the system gains usable information about the feasibility of using weather-affected fixes.

We next formalize these roles by introducing utility models for each stakeholder, which capture trade-offs between delay reduction, operational risk, and airline-specific costs in pathfinder acceptance decisions. First, we define a flight model that captures how individual flights decide whether to accept a pathfinder offer. We then specify an ATC utility that reflects system objectives such as reopening a closed departure fix and reducing overall delay. Finally, we describe an airline dispatcher utility that encodes carrier-specific priorities, including delay costs and operational regulations. These utility models provide a common framework to analyze how stakeholders make decisions, how pathfinder candidates are selected, and how offer sequences are formed in practice.

\subsubsection{Flight's Utility Model}
\label{sec:flight_model}
To describe pathfinder acceptance at the individual flight level, we formulate a stylized, utility-based decision model. We intentionally simplified the models to preserve analytical tractability, focusing on the key trade-offs that shape flight acceptance behavior. 

We consider a set of candidate flights $N = \{1, 2, \dots, n\}$, where each flight $i \in N$ independently determines whether to participate in the pathfinder operation. No coordination or information sharing among flights is assumed. Each candidate flight faces a binary participation choice given by
\begin{equation}
x_i =
\begin{cases}
1, & \text{if flight } i \text{ accepts pathfinding offer,} \\
0, & \text{otherwise.}
\end{cases}
\end{equation}

Each flight evaluates whether to accept a pathfinder offer by weighing three components: a participation cost, a failure cost, and a reward. The participation cost $c_i(x_i)$ represents the operational or perceived burden associated with undertaking the pathfinder attempt, with $c_i(0) = 0$ when the flight declines. Examples include increased tactical weather-avoidance workload and handling of nonstandard clearances. The failure cost $d_i(x_i)$ captures the expected loss incurred if a pathfinding attempt does not succeed and likewise satisfies $d_i(0) = 0$. In this context, a failure corresponds to the situation in which a flight accepts the assignment and departs but is unable to trigger fix reopening because adverse weather conditions persist. Such unsuccessful attempts may lead to additional fuel consumption, route deviations, and passenger discomfort due to turbulence, as reported by aviation stakeholders during our interviews \cite{interview}. The reward $T_i(x_i)$ is received whenever the flight accepts the offer and departs, independent of the eventual success of the pathfinder operation. This term represents the flight’s incentive to avoid further ground delay by departing promptly, and is interpreted as a form of incentivization for assuming the pathfinder attempt.

The flight’s decision determines whether a pathfinding attempt is undertaken. If flight $i$ accepts the pathfinder offer ($x_i = 1$), the attempt succeeds with probability $P_{\mathrm{success},\,i}$; if it declines ($x_i = 0$), no pathfinding attempt is made and the probability of success is zero. Accordingly, the probability of success can be written as
\begin{equation}
P_{\mathrm{success},\,i}(x_i) =
\begin{cases}
P_{\mathrm{success},\,i}, & \text{if } x_i = 1, \\
0, & \text{if } x_i = 0.
\end{cases}
\end{equation}
The resulting utility for flight \( i \) is expressed as
\begin{equation}
U_i(x_i) =
\begin{cases}
T_i - c_i - (1 - P_{\mathrm{success},\,i})d_i, & \text{if } x_i = 1, \\
0, & \text{if } x_i = 0.
\end{cases}\label{eq:utility}
\end{equation}
The utility function in \eqref{eq:utility} captures the flight’s overall expected payoff from its participation decision. When the flight accepts the pathfinder assignment ($x_i = 1$), it obtains a reward $T_i$, incurs a participation cost $c_i$, and is exposed to an expected loss from potential failure weighted by the mission success probability. When the flight chooses not to participate ($x_i = 0$), the resulting utility is zero.

\begin{figure}[t]
    \centering
    \includegraphics[width=0.8\linewidth]{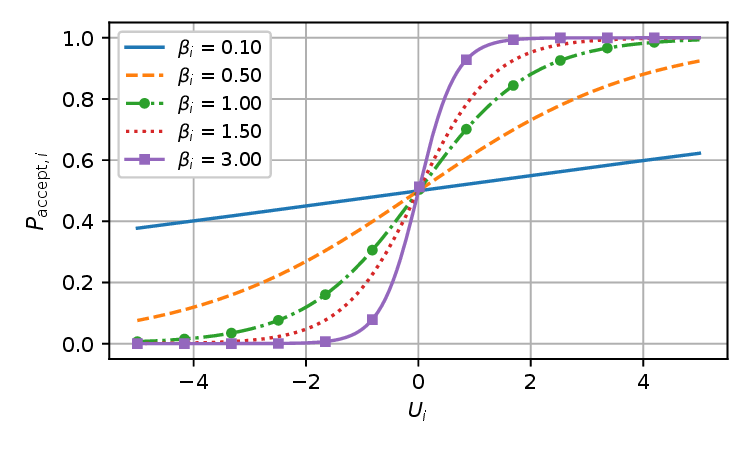}
    \caption{\(P_{\mathrm{accept},\,i}\) as a function of utility \( U_i \) for different sensitivity values \( \beta_i \). Higher \( \beta_i \) yields more deterministic decisions, while lower \( \beta_i \) leads to more randomized acceptance.}
    \label{fig:flight_utility}
\end{figure}

The acceptance probability is represented using a logistic function, which is widely adopted in discrete choice modeling to account for bounded rationality through stochastic utility maximization \cite{logit_model}:
\begin{equation}\label{eq:accept_prob}
P_{\mathrm{accept},\,i}(x_i) = \frac{1}{1 + e^{-\beta_i U_i(x_i)}},
\end{equation}
where \( \beta_i > 0 \) controls the sensitivity of the flight’s choice to utility differences, with larger values yielding more deterministic behavior and smaller values introducing greater randomness. This effect is visualized in \Cref{fig:flight_utility}.

\subsubsection{ATC's Utility Model}
\label{sec:atc_model}
From the ATC perspective, pathfinder selection is a strategic decision to minimize overall system delay and the risk of unsuccessful or rejected pathfinder attempts. The controller selects a candidate flight based on its expected contribution to reducing system-wide ground delay. To capture this, we define an ATC's utility for each candidate flight $i \in N$ as
\begin{equation}
U^{\mathrm{ATC}}_i
=
P_{\mathrm{accept},i} \cdot \bigg(
\Delta D^{\mathrm{sys}}_i - \lambda_\mathrm{ATC} \, G^\mathrm{ATC}_i
\bigg),
\label{eq:atc_model}
\end{equation}
where the term \(P_{\mathrm{accept},i}\) is the acceptance probability of flight \(i\), defined in \Cref{eq:accept_prob}, and reflects how likely a flight is to accept a pathfinder offer from ATC. The quantity \(\Delta D_{\mathrm{sys},i} \ge 0\) represents the expected system delay reduction if flight \(i\) serves as the pathfinder. This term can depend on factors such as the flight's position in the departure queue or the number of downstream flights that could benefit from reopening the fix.

The term \(G^{\mathrm{ATC}}_i \ge 0\) captures the operational burden or penalty, from the ATC viewpoint, associated with selecting flight \(i\) as the pathfinder. Concrete examples of operational burden from the ATC perspective include: additional coordination efforts between stakeholders (e.g., ATC communicating with air carriers), departure sequence disruption (e.g., reshuffling aircraft on airport taxiways), or added sector complexity (e.g., from the potential of flights to deviate due to worse-than-expected weather) The weighting parameter \(\lambda_\mathrm{ATC} \ge 0\) controls how strongly these penalties are traded off against potential delay savings. As $\lambda_{\mathrm{ATC}}$ increases, the ATC objective increasingly prioritizes reducing operational burden over delay reduction,  favoring candidates that are easier to accommodate operationally even if they yield somewhat smaller delay reductions.

Putting these terms together, \(U^{\mathrm{ATC}}_i\) can be interpreted as the expected net gain of selecting flight \(i\) as the pathfinder from the ATC perspective. Flights with high acceptance probability, large expected delay reduction, and low operational penalty obtain higher utility values. 

\subsubsection{Airline Dispatcher's Utility Model}\label{sec:dispatcher_model}
From the airline dispatcher's perspective, the main objective is to reduce carrier specific costs while avoiding operational risk. Typical concerns include DOT3 rules, missed passenger connections, crew duty and legality constraints, and disruptions to aircraft rotations. To capture this, we define the dispatcher-side utility for each candidate flight \(i\) belonging to a given airline. 

Let \({A \subseteq N}\) denote the set of flights from a specific airline that are eligible to be considered for pathfinding. For flight \(i \in {A}\), the dispatcher utility used for sequencing pathfinder offer is
\begin{equation}
U^{\mathrm{disp}}_i
=
P_{\mathrm{accept},i} \cdot \bigg(
B^{\mathrm{dep}}_i - \lambda_\mathrm{disp} \, G^\mathrm{disp}_i
\bigg),
\qquad i \in {A}.
\label{eq:disp_model}
\end{equation}
The term \(B^{\mathrm{dep}}_i \ge 0\) represents the expected departure schedule benefit to the airline if flight \(i\) serves as the pathfinder. It includes, for example: avoiding a potential DOT3 violation, reducing expected downstream delay for that aircraft, preserving crew schedules, and preserving passenger connections on later legs. Note that we consider detailed models (e.g., large-scale integer programs \cite{Hane1995FleetAssignment}) of airline tail recovery, crew scheduling, and related problems to be out-of-scope for this study: Such microscopic models may be incorporated as inputs which inform $B^\mathrm{dep}_i$. This is particularly relevant because unsuccessful or prolonged pathfinding attempts may keep an aircraft waiting on the ground for an extended period, increasing the risk of exceeding the three-hour tarmac delay limit. 

The term \(G^{\mathrm{disp}}_i \ge 0\) captures dispatcher perceived risk or cost associated with selecting flight \(i\) as the pathfinder. Examples include insufficient turnaround readiness, uncertainty about maintenance or fueling status, tight crew duty limits, or a high risk of causing additional delay on subsequent legs. As noted above, these factors are represented at an aggregate level rather than through explicit operational models. The parameter \(\lambda_{\mathrm{disp}} \ge 0\) controls how strongly these risks are penalized relative to the expected benefit. Larger values of $\lambda_{\mathrm{disp}}$ increase the emphasis on dispatcher-side risk, favoring flights that are easier to accommodate operationally even if they provide smaller departure benefits.

Under this formulation, \(U^{\mathrm{disp}}_i\) represents how beneficial it is for the airline to select flight \(i\) as a pathfinder candidate. Flights with higher acceptance probability, larger departure side benefit, and lower operational risk receive higher utility values. 

\subsection{Worst Case Scenario Analysis}
\label{sec:worst_case}
Having defined the utility structures for all relevant stakeholders and described how their decision rules shape the pathfinder selection process, we now turn to the resulting implications at the system level. From an operational standpoint, one of the most disruptive outcomes arises when no flight agrees to assume the pathfinder attempt, leaving the system without the information required to reopen a weather-impacted departure fix. While pathfinder attempts can also fail to be initiated due to ATC or airline-side constraints, our worst-case analysis here focuses on the risk arising from collective rejection by candidate flights. This motivates a focused examination of a critical scenario in which \emph{all} candidate flights decline the pathfinder assignment, which is operationally severe because it prevents fix reopening and propagates departure delays.

To analyze this failure mechanism, we consider a sequential decision framework in which individual flights are approached one by one with a pathfinder offer. Emphasis is placed on the worst case outcome, where every flight rejects the offer, leading to a complete breakdown of the pathfinding process. Each flight’s response is modeled probabilistically. The objective is to characterize how the likelihood of this collective rejection event varies under different system assumptions. We examine a baseline configuration along with two extended models that incorporate selfless decision-making and environmental uncertainty. For analytical tractability, simplifying assumptions are introduced and stated explicitly within each subsection.

\subsubsection{Baseline Failure Mode}
\label{sec:independent_assumption}
We first consider a baseline configuration in which each flight makes an independent accept-or-reject decision based solely on its own utility. Flights do not observe the acceptance decisions of other flights and make their choices independently according to the logistic acceptance model introduced in \Cref{sec:flight_model}, given by \Cref{eq:accept_prob}. The corresponding rejection probability is given by
\begin{equation}
P_{\text{reject},\,i}(x_i) = 1 - P_{\text{accept},\,i}(x_i) = \frac{1}{1 + e^{\beta_i U_i(x_i)}}.
\end{equation}

We partition the flight population into two behavioral types, labeled as \emph{rejective} and \emph{receptive}. Each flight is independently assigned to the rejective class with probability $\alpha$, and to the receptive class with probability $1 - \alpha$. To facilitate analysis at the aggregate level, these two groups are characterized by representative utility values $U^- < 0$ and $U^+ > 0$, respectively. The corresponding rejection probabilities for the two groups are given by
\begin{equation}
P_{\text{rejective}} = \frac{1}{1 + e^{\beta U^-}}, \quad
P_{\text{receptive}} = \frac{1}{1 + e^{\beta U^+}},
\end{equation}
with $P_{\text{rejective}} > 0.5$ and $P_{\text{receptive}} < 0.5$. For simplicity, a common sensitivity parameter $\beta$ is assumed for all flights. This formulation results in a heterogeneous population with two flight types: those that tend to accept pathfinder offers and those that tend to reject them.

The failure event corresponds to the case in which all $n$ flights decline the pathfinder offer. Let $W(\alpha)$ denote the probability of this outcome, where $\alpha \in [0,1]$ denotes the proportion of rejective flights in the overall population. The probability is given by
\begin{equation}
W(\alpha) = \sum_{k=0}^{n} \binom{n}{k} \alpha^k (1 - \alpha)^{n - k}
(P_{\text{rejective}})^k (P_{\text{receptive}})^{n - k}.
\end{equation}
Using the binomial theorem, this simplifies to
\begin{equation}
W(\alpha) = \left(\alpha P_{\text{rejective}} + (1 - \alpha) P_{\text{receptive}}\right)^n.
\end{equation}
Since $P_{\text{rejective}} > P_{\text{receptive}}$, the failure probability $W(\alpha)$ is an increasing function of $\alpha$, indicating that system vulnerability grows as the share of rejective flights increases. The sensitivity of $W(\alpha)$ to changes in $\alpha$ becomes stronger when the difference between utilities $U^-$ and $U^+$ increases (i.e., $U^-\ll U^+$) and when the sensitivity parameter $\beta$ increases, whereas larger values of $n$ reduce the overall likelihood of complete rejection.

We introduce a critical threshold $\alpha^*$ defined implicitly by the condition $W(\alpha^*) = \delta$, where $\delta \in (0,1)$ denotes the maximum tolerable probability of complete rejection. Solving this equation gives
\begin{equation}
\alpha^* = \frac{\delta^{1/n} - P_{\text{receptive}}}{P_{\text{rejective}} - P_{\text{receptive}}}.
\end{equation}
This threshold separates stable and fragile conditions. When $\alpha \geq \alpha^*$, the probability of complete rejection becomes substantial, indicating system failure. Operationally, crossing this threshold means that most available flights are unlikely to accept the pathfinder offer. As a result, reopening attempts are more likely to fail, and departure capacity is more likely to remain unavailable for extended periods during bad weather.

\subsubsection{Considering Flight's Selfless Behavior}
\label{sec:selfless_behavior}
We next extend the baseline model to account for the possibility of selfless behavior in pathfinder decisions. In operational settings, not all flights assess pathfinder offers purely from a self-interested perspective. Some flights may also weigh broader operational benefits of acceptance, such as facilitating departures for other delayed flights. This extension introduces a limited form of interdependence by allowing flights to account for system benefits in their acceptance decisions, while preserving the independent decision structure of the baseline model.

To model this effect, we introduce a selfishness parameter $S_i$ for each flight $i$, with $S_i \in [0,1]$. A value of $S_i = 1$ represents a fully selfish flight that considers only its own utility, whereas $S_i = 0$ corresponds to a fully selfless flight that fully internalizes the system outcomes of its decision.

The adjusted utility for flight $i$ is constructed by augmenting the individual utility with a system contribution that accounts for the anticipated effect of collective rejection,
$U_i^{\text{sys}} = U_i(x_i) + (1 - S_i)\,\gamma_i R$, where $\gamma_i > 0$ denotes a sensitivity coefficient and $R$ represents a perceived system rejection risk used by each flight in its decision-making. While $R$ is not an explicit argument of $W(\alpha)$, it is assumed to be positively associated with the system-wide probability of complete rejection characterized by $W(\alpha)$. This additional term captures the flight’s awareness of the risk associated with widespread rejection and embeds a notion of responsibility for system-wide improvements into the decision-making process.

Based on \( U_i^{\text{sys}} \), the rejection probability becomes 
\begin{equation} 
    P_{\text{reject}}(x_i) = \frac{1}{1 + e^{\beta_i U_i^{\text{sys}}}}. 
\end{equation}
Under the extension that incorporates selfless behavior, the probability of the worst case event is given by
\begin{equation}
W_{sys}(\alpha) = \left( \alpha P_{\text{rejective}}^{sys} + (1 - \alpha) P_{\text{receptive}}^{sys} \right)^n,
\end{equation}
where
\begin{align}
P_{\text{rejective}}^{sys} &= \frac{1}{1 + e^{\beta (U^- + (1 - S)\gamma R)}}, \text{ and} \\
P_{\text{receptive}}^{sys} &= \frac{1}{1 + e^{\beta (U^+ + (1 - S)\gamma R)}}.
\end{align}
For analytical convenience, all flights are assumed to share a common selfishness level $S$ and a common system-influence parameter $\gamma$. Lower values of $S$ corresponds to greater selflessness, and larger values of $\gamma$ reflects stronger sensitivity to system-wide failure, both acting to decrease rejection probabilities and enhance overall system resilience.

\subsubsection{Considering Uncertainty in Environments}
\label{sec:environment_uncertainty}
Flight decisions are also subject to external sources of uncertainty, including imperfect information (e.g., specific dispatcher-pilot interactions) and random disturbances (e.g., in terms of coordination lag times). To represent this effect, we augment the utility function with an additive stochastic noise term $\xi$, modeled as a zero-mean random variable with noise intensity $\theta$, which is assumed to be common across all flights. This formulation reflects the assumption of a shared realization of external uncertainty affecting all decision-makers simultaneously. Two alternative noise specifications are considered, namely Gaussian and Rademacher distributions, enabling the analysis of both continuous and discrete forms of randomness and their influence on system behavior. Gaussian noise can be viewed as a stylized representation of continuous uncertainty in weather information, such as imprecision in estimating the location or extent of convective weather. Rademacher noise represents a coarse binary perturbation that reflects threshold-based operational assessments of whether weather impacts are present.

Environmental uncertainty is incorporated into the flight’s decision-making by the modified utility $U_i^{\text{env}} = U_i(x_i) + \xi, 
$ where $\xi$ denotes a zero-mean random variable with scale parameter $\theta \in \{\sigma, \kappa\}$. Specifically, $\theta = \sigma$ when $\xi \sim \mathcal{N}(0,\sigma^2)$, and $\theta = \kappa$ when $\xi \sim \text{Rademacher}(\pm \kappa)$. The unified parameter $\theta$ is used to provide a consistent representation across both noise distributions.

Under environmental uncertainty, the worst case probability becomes a function of both the rejective flight proportion $\alpha$ and the noise scale $\theta$, which we denote by $W(\alpha,\theta)$. The corresponding tipping point $\alpha^*(\theta)$ is defined implicitly as the value of $\alpha$ satisfying $W(\alpha^*(\theta), \theta) = \delta,$ for a prescribed threshold $\delta \in (0,1)$.

To examine how the tipping point $\alpha^*(\theta)$ varies with respect to the noise scale $\theta$, we invoke the implicit function theorem by defining the auxiliary function
\begin{equation}
G(\alpha(\theta), \theta) := W(\alpha(\theta), \theta) - \delta .
\end{equation}
Under this formulation, the tipping point is characterized by the root condition $G(\alpha^*(\theta), \theta) = 0$. This representation enables differentiation of $\alpha^*(\theta)$ with respect to $\theta$, treating the tipping point as an implicitly defined function. The resulting derivative is
\begin{equation}
\frac{\mathrm{d}\alpha^*}{\mathrm{d}\theta}
= - \frac{\partial G / \partial \theta}{\partial G / \partial \alpha}
= - \frac{\partial W / \partial \theta}{\partial W / \partial \alpha}.
\end{equation}
This expression follows directly from application of the implicit function theorem to the relation $G(\alpha^*(\theta), \theta) = 0$.

Since $W(\alpha,\theta)$ is an increasing function of $\alpha$, a higher proportion of rejective flights necessarily raises the likelihood of complete rejection, which implies $\partial W / \partial \alpha > 0$. Consequently, the sign of $\mathrm{d}\alpha^* / \mathrm{d}\theta$ is determined entirely by the sign of $\partial W / \partial \theta$, which quantifies the influence of environmental uncertainty on the worst case probability. The derivative $\mathrm{d}\alpha^* / \mathrm{d}\theta$ therefore indicates how the tipping point $\alpha^*$ shifts in response to changes in uncertainty. When $\partial W / \partial \theta > 0$, increasing uncertainty amplifies the failure probability, leading to a decrease in $\alpha^*$ and rendering the system more fragile. Conversely, when $\partial W / \partial \theta < 0$, greater uncertainty mitigates the risk of complete rejection, resulting in an increase in $\alpha^*$ and a corresponding enhancement in system robustness.

The gradient-based analysis above characterizes how environmental uncertainty shifts the system’s tipping point between stable and fragile regimes. In \Cref{sec:worst_result}, we show that operational settings can fall on either side of this boundary, implying that system performance can be sensitive to which flights accept pathfinder offers.  This sensitivity motivates the need for decision support mechanisms that actively shape pathfinder acceptance outcomes. In particular, since acceptance likelihood depends on the offer order, the sequencing of pathfinder offers becomes an important operational control.

\subsection{Sequence Optimization in Pathfinder Selection}\label{sec:sequence_optimization}
We now turn to the question of how pathfinder offers should be allocated in practice. Instead of relying on a predetermined flight schedule, we construct offer sequences that use stakeholder utilities to guide the ordering while also accounting for the risk that all candidate flights may reject the pathfinder offer. We formulate optimization problems that generate candidate sequences under both ATC-initiated and dispatcher-initiated pathfinder operations. To provide the optimization with realistic parameter values (e.g., expected delay savings and operational penalties), we develop a discrete-event simulation of the airport departure process for the case study airport. The simulation evaluates the operational impact of assigning different flights as pathfinders and links the utility models to sequencing decisions.

\subsubsection{Optimization Formulations}\label{sec:formulations}
\paragraph{ATC-Initiated Formulation}
We first present an optimization formulation for ATC-initiated pathfinder selection, in which the ATC constructs the offer sequence and determines which flights are approached to be a pathfinder and in what order. The model selects and orders candidate flights to maximize the expected system benefit of the pathfinder process while accounting for acceptance behavior, operational costs, and coordination constraints. 
\begin{align}
\max_{x,y} \quad &
\sum_{k=1}^n \sum_{i=1}^n x_{i,k}\, p_{i,k}
\Big( \Delta D^{\text{sys}}_{i,k} - \lambda_\mathrm{ATC} G^{\mathrm{ATC}}_{i,k} \Big)
\,\cdot \notag\\
&\prod_{j=1}^{k-1} \left(1 - \sum_{i'=1}^n x_{i',j} p_{i',j}\right),
\label{eq:atc_obj} \\
\text{s.t.} \quad
& \sum_{i=1}^n x_{i,k} = y_k, \quad \forall k=1,\dots,n,
\label{eq:atc_assign} \\
& y_k \ge y_{k+1}, \quad \forall k=1,\dots,n-1,
\label{eq:atc_prefix} \\
& \sum_{k=1}^n x_{i,k} \le 1, \quad \forall i=1,\dots,n,
\label{eq:atc_unique} \\
& \sum_{i=1}^n \sum_{k=1}^n e_{i,k}\, x_{i,k} \le B,
\label{eq:atc_budget} \\
& x_{i,k} \in \{0,1\},\; y_k \in \{0,1\},  \quad \forall i,k.
\label{eq:atc_binary}
\end{align}

For notational clarity, we denote by \(p_{i,k}\) the acceptance probability of flight \(i\) when offered at position \(k\), extending the flight acceptance model introduced earlier. Here, position \(k\) refers to the order in which a flight is offered the pathfinder attempt within a given offer sequence. In this formulation, all flights in the system \(N\) are considered potential candidates. For each flight \(i\) and sequence position \(k\), the binary variable \(x_{i,k}\) indicates whether flight \(i\) is assigned to position \(k\), while \(y_k\) specifies whether the \(k\)-th position in the sequence is active. 

The objective \Cref{eq:atc_obj} maximizes the expected net system benefit of the entire offer sequence. The term \(p_{i,k}(\Delta D^{\mathrm{sys}}_{i,k} - \lambda_{\mathrm{ATC}}G^{\mathrm{ATC}}_{i,k})\) corresponds directly to the ATC utility model, in which system delay reduction is traded off against operational or readiness costs. If flight \(i\) is placed at position \(k\), this expected contribution is weighted by the probability that the process reaches position \(k\), captured by the product term since any earlier acceptance stops the sequence.

The constraints enforce a valid sequencing structure. \Cref{eq:atc_assign} requires that each active position contain exactly one flight, while \Cref{eq:atc_prefix} ensures that positions are used consecutively starting from the first position (i.e., if position \(k\) is used, then all positions \(1,\ldots,k-1\) must also be used). \Cref{eq:atc_unique} guarantees that each flight can appear at most once in the sequence. The coordination budget in \Cref{eq:atc_budget} limits total communication or workload cost associated with issuing offers. Finally, \Cref{eq:atc_binary} enforces the discrete nature of the assignment decisions. Overall, this formulation yields an ATC-prioritized offer sequence that balances acceptance likelihood, system delay reduction, operational cost, and coordination constraints.

\paragraph{Dispatcher-Initiated Formulation}
In the dispatcher-initiated formulation, the sequencing problem is defined over the airline-specific set of flights \(A\). The binary variable \(x_{i,k}\) indicates whether flight \(i \in A\) is assigned to position \(k\), while \(y_k\) denotes whether position \(k\) is active. The dispatcher constructs an offer sequence over this set by solving the optimization problem.
\begin{align}
\max_{x,y} \quad &
\sum_{k=1}^{|A|} \sum_{i \in A} x_{i,k}\, p_{i,k}
\Big( B^{\text{dep}}_{i,k} - \lambda_\mathrm{disp} G^{\mathrm{disp}}_{i,k} \Big)
\,\cdot \notag\\
&\prod_{j=1}^{k-1} \left(1 - \sum_{i' \in A} x_{i',j} p_{i',j}\right),
\label{eq:disp_obj}
\\
\text{s.t.} \quad 
& \sum_{i \in A} x_{i,k} = y_k,
\quad \forall k = 1,\dots,|A|,
\label{eq:disp_assign}
\\
& y_k \ge y_{k+1},
\quad \forall k = 1,\dots,|A|-1,
\label{eq:disp_prefix}
\\
& \sum_{k=1}^{|A|} x_{i,k} \le 1,
\quad \forall i \in A,
\label{eq:disp_unique}
\\
& \sum_{i \in A} \sum_{k=1}^{|A|} e_{i,k}\, x_{i,k} \le B,
\label{eq:disp_budget}
\\
& x_{i,k} \in \{0,1\}, \; y_k \in \{0,1\},
\quad \forall i,\; k.
\label{eq:disp_binary}
\end{align}

The objective \Cref{eq:disp_obj} captures airline-specific priorities: placing flight \(i\) at position \(k\) yields an expected contribution \(p_{i,k}(B^{\mathrm{dep}}_{i,k} - \lambda_{\mathrm{disp}} G^{\mathrm{disp}}_{i,k})\), where \(B^{\mathrm{dep}}_{i,k}\) reflects departure-side operational benefits and \(G^{\mathrm{disp}}_{i,k}\) represents dispatcher-perceived risks such as readiness concerns or downstream schedule disruption. This expression mirrors the dispatcher utility model introduced earlier, in which benefits and operational risks are combined into a single decision metric. As in the ATC-initiation formulation, the contribution is weighted by the probability of reaching position \(k\), since any earlier acceptance terminates the process. 





The dispatcher-initiated formulation uses the same structural constraints as the ATC-initiated version, including prefix, uniqueness, budget, and integrality constraints. The key difference lies in the objective function, which prioritizes airline operational benefits rather than system delay reduction, producing a dispatcher-prioritized pathfinder sequence.


\subsubsection{Discrete Event Simulation for Parameter Calculation}\label{sec:des_calculation}

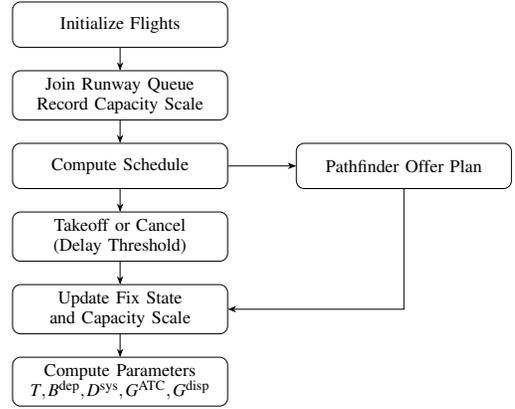
\begin{figure}[t]
\centering
\resizebox{0.75\linewidth}{!}{%
\begin{tikzpicture}[
  font=\small,
  node distance=4mm and 1.5mm,
  box/.style={rectangle, rounded corners, draw, align=center, minimum width=3.8cm, minimum height=8mm},
  >={Stealth[length=1.4mm]}
]

\node[box] (gen) {Initialize Flights};
\node[box, below=of gen] (queue) {Join Runway Queue\\Record Capacity Scale};
\node[box, below=of queue] (sched) {Compute Schedule};
\node[box, below=of sched] (takeoff) {Takeoff or Cancel\\(Delay Threshold)};
\node[box, below=of takeoff] (fix) {Update Fix State\\and Capacity Scale};
\node[box, below=of fix] (metrics) {Compute Parameters\\\(T, B^{\text{dep}}, D^{\text{sys}}, G^\mathrm{ATC},G^\mathrm{disp}\)};

\node[box, right=12mm of sched] (offer) {Pathfinder Offer Plan};

\draw[->] (gen) -- (queue);
\draw[->] (queue) -- (sched);
\draw[->] (sched) -- (takeoff);
\draw[->] (takeoff) -- (fix);
\draw[->] (fix) -- (metrics);
\draw[->] (sched) -- (offer);
\draw[->] (offer) |- (fix);
\end{tikzpicture}}
\caption{Simulation flow for runway sequencing, fix updates, and pathfinder evaluation.}
\label{fig:sim_flow}
\end{figure}

To obtain realistic values for five parameters \((T,\; B^{\text{dep}},\; D^{\text{sys}},\; G^\mathrm{ATC}, \text{ and }G^\mathrm{disp})\) used in the utility models and sequencing formulations, we develop a discrete-event simulation of departure operations at a major hub airport (e.g., JFK). The simulator is implemented in Python using SimPy and models flights, runways, and departure fixes as interacting processes evolving in continuous time. Time is measured in minutes from \(t=0\), and all separation and policy rules are defined on this scale. A high-level overview of the simulation workflow is shown in \Cref{fig:sim_flow}.

Each flight is represented with operational attributes (e.g., scheduled push time, unimpeded taxi time, wake category, destination region) and simulation-state attributes (e.g., pathfinder flag, cancellation status, connection value, crew duty limit). During a run, a flight waits until it is ready, undergoes a stochastic taxi process, selects a departure fix and runway, joins the corresponding queue, and ultimately departs or cancels if a delay threshold is exceeded. The simulator records join time, takeoff time, assigned runway and fix, and relevant delay and benefit metrics.

Runway operations are modeled as queues that release flights according to wake-turbulence spacing and the current departure capacity. For a leading wake class $\ell$ and trailing class $\tau$, the base headway is
\begin{equation}
h_{\text{base}}(\ell,\tau)=\frac{s_{\ell,\tau}+r}{60},
\end{equation}
where \(s_{\ell,\tau}\) is the required wake separation in seconds and \(r\) is an additional roll buffer. This base spacing is scaled by a factor \(sigma\) that depends on how many departure fixes are open, giving an effective headway \(h_{\text{eff}} = h_{\text{base}} \sigma\). Opening more fixes reduces $\sigma$ and increases runway throughput. A minimum spacing can be enforced for pathfinder flights. During the simulation, the runway updates predicted departure times using these headways, advances time, tracks the most recent wake class, and cancels flights whose delay exceeds a threshold.

Departure fixes are treated as binary open/closed resources managed by the simulation logic. Whenever a fix opens, the simulation updates the capacity scale based on the number of available fixes and records the event. Flights choose a fix using a simple priority rule: they first try a destination-preferred fix, then a runway-preferred fix, and if neither is available, they select among the open fixes in a round-robin manner. Flights may still join a runway queue when all fixes are closed, but they cannot depart until at least one fix becomes usable. 

The simulator also includes a pathfinder module for evaluating policies that allow queue jumping and fix opening. A pathfinder plan specifies which flight may become the pathfinder, when the offer sequence begins, how much overhead to apply for declines and acceptance, which fixes should open after a successful pathfinder departure, and what event triggers that opening (e.g., the pathfinder’s takeoff). When the plan is active, the simulator runs the offer sequence, designates the pathfinder if it is still eligible, reinserts it into the runway queue at the position yielding the earliest predicted departure, and opens the specified fixes when the trigger occurs. The capacity scale is then updated to reflect the new fix configuration.

To measure policy impact, the simulator performs paired baseline and pathfinder runs using a common random seed. Let \(w_f\) denote the waiting time of flight \(f\), limited to a maximum value. The system delay improvement is \(\Delta D^{\text{sys}} = \sum_f w^{\text{base}}_f - \sum_f w^{\text{pf}}_f.\) For the designated pathfinder, the simulator reports its flight reward \(T\), connection gain, and airline benefit \(B^{\text{dep}}\). Coordination costs are also computed: an ATC cost proportional to the number of positions the pathfinder jumps and a dispatcher cost determined by the flights it overtakes. These outputs provide the quantities required by the sequence optimization formulations, and in practical deployments would be calibrated using operational data rather than simulation parameters.

\section{Numerical Results}\label{sec:results}
Building on the formal decision-making models developed earlier in \Cref{sec:modeling_markov_chain,sec:pathfinder_selection,sec:worst_case,sec:sequence_optimization}, we now present results that illustrate how the proposed models can be used to analyze operational behavior. We use the analytical framework and simulation tools developed earlier to examine system performance, vulnerability, and pathfinder offer sequence decisions in realistic settings.

\subsection{Markov Chain–Based System Behavior}\label{sec:markov_result}
\subsubsection{Steady-State Behavior and Long-Term Dynamics}
\label{sec:steady_state}
\begin{figure*}[htbp]
  \centering
  \includegraphics[width=\linewidth]{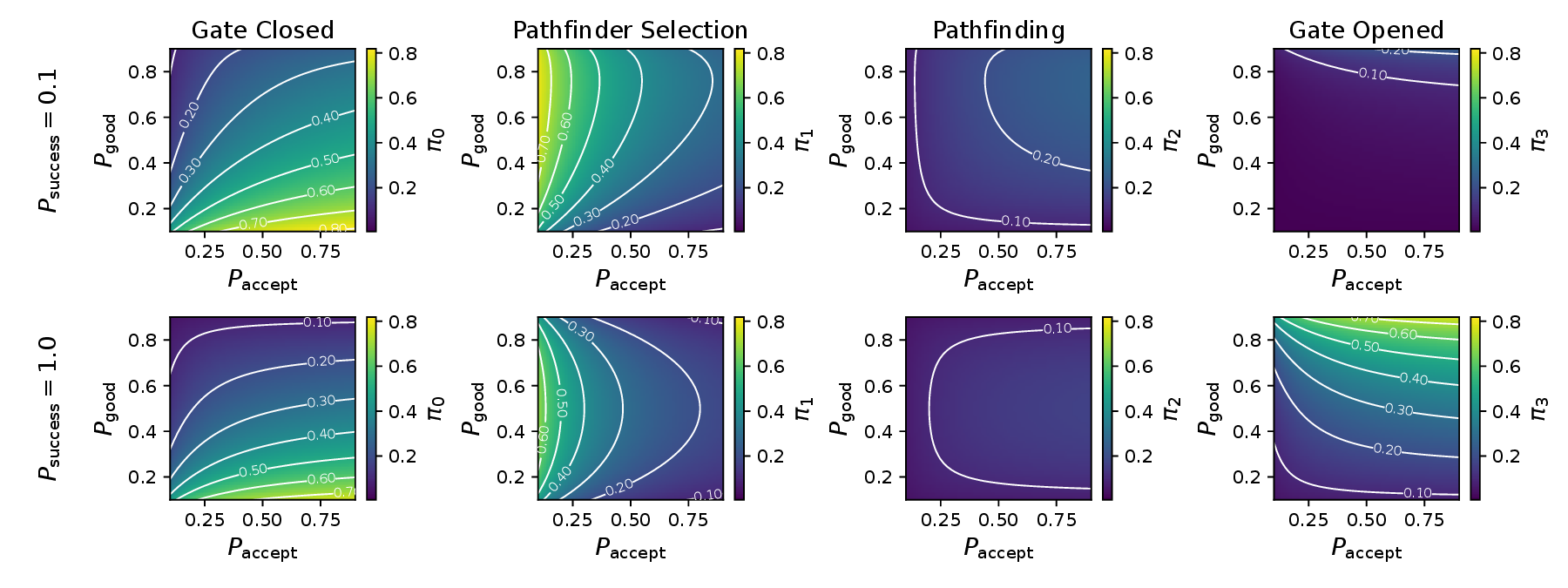}
  \caption{Steady-state distribution (\(\pi_0\) to \(\pi_3\)) of a four-state Markov chain under varying values of \(P_{\mathrm{good}}\) and \(P_{\mathrm{accept}}\), for two different values of \(P_{\mathrm{success}}\): 0.1 (top row) and 1.0 (bottom row). Each column corresponds to one of the four system states: \textit{Gate Closed}, \textit{Pathfinder Selection}, \textit{Pathfinding}, and \textit{Gate Opened}. The color bar represents the steady-state probability \(\pi_i\) for each state. The same color scale is used within each column for fair comparison between the two \(P_{\mathrm{success}}\) settings.}
  \label{fig:steady_state}
\end{figure*}

In a Markov chain, the steady state describes how the system behaves in the long run, meaning the proportion of time it spends in each state once short-term fluctuations fade. This provides a clear picture of which operational conditions dominate over repeated pathfinding cycles. We compute the distribution numerically in our analysis. The sensitivity analysis results are shown in \Cref{fig:steady_state}, illustrating how long-run occupancy varies across the four system states under different values of \(P_{\mathrm{good}}\), \(P_{\mathrm{accept}}\), and \(P_{\mathrm{success}}\). \Cref{fig:steady_state} compares scenarios with low and high pathfinding success, highlighting how input parameters affect steady-state probabilities.

As \(P_{\mathrm{good}}\) increases, the system tends to remain in the \textit{Gate Opened} state (\(\pi_3\)), while the \textit{Gate Closed} state (\(\pi_0\)) becomes less frequent. Increasing \(P_{\mathrm{accept}}\) leads to more transitions into pathfinding, reflected in lower \(\pi_1\) and \(\pi_2\). A key relationship, \(\pi_2 = P_{\mathrm{accept}} \pi_1\), ensures that pathfinding (State~2) always has a lower steady-state probability than \textit{Pathfinder Selection} (State~1), explaining why the system spends relatively little time in State~2. This effect is particularly pronounced when \(P_{\mathrm{success}} = 1\), in which case the system primarily alternates between the \textit{Gate Closed} and \textit{Gate Opened} states. When pathfinding success is low (\(P_{\mathrm{success}} = 0.1\)), the system becomes dominated by \textit{Gate Closed} and selection states, while the \textit{Gate Opened} state (\(\pi_3\)) sharply decreases. The rise in \(\pi_2\) reflects its dependence on \(\pi_1\), and maintaining throughput requires high values of both \(P_{\mathrm{good}}\) and \(P_{\mathrm{accept}}\).

This steady-state analysis highlights that the system's ability to sustain nominal operations depends not only on environmental conditions but also on decision-making parameters, particularly \( P_{\mathrm{accept}}\), which determines how readily the system transitions into pathfinding. As \( P_{\mathrm{accept}}\) increases, the system is more responsive to opening fixes, reducing time spent in intermediate states. These results underscore the importance of explicitly modeling pathfinder selection decisions, validating our focus on decision-level mechanisms that shape \(P_{\mathrm{accept}}\).

\subsubsection{Operational Performance Implications}\label{sec:markov_operation}
\begin{figure*}[t]
    \centering
    \begin{subfigure}{0.28\linewidth}
        \centering
        \includegraphics[width=\linewidth]{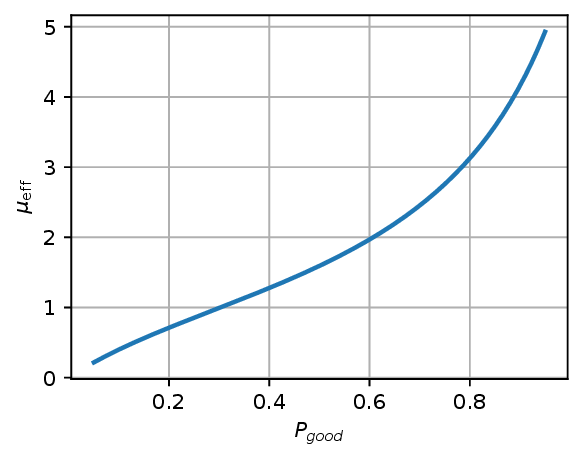}
        \caption{$\mu_{\mathrm{eff}}$ vs.\ $P_{\mathrm{good}}$}
        \label{fig:markov_operation_a}
    \end{subfigure}%
    \hfill%
    \begin{subfigure}{0.28\linewidth}
        \centering
        \includegraphics[width=\linewidth]{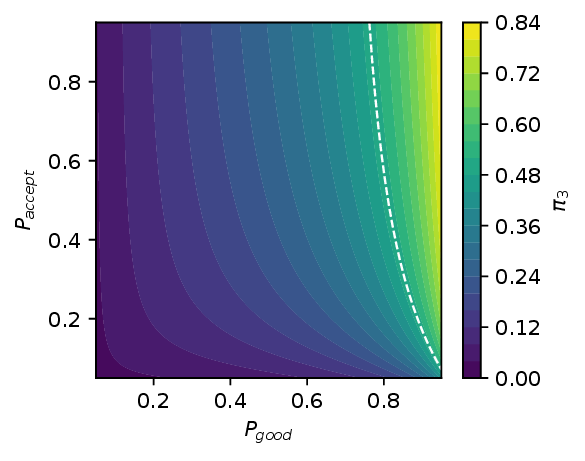}
        \caption{$\pi_{3}$ over $(P_{\mathrm{good}}, P_{\mathrm{accept}})$}
        \label{fig:markov_operation_b}
    \end{subfigure}%
    \hfill%
    \begin{subfigure}{0.28\linewidth}
        \centering
        \includegraphics[width=\linewidth]{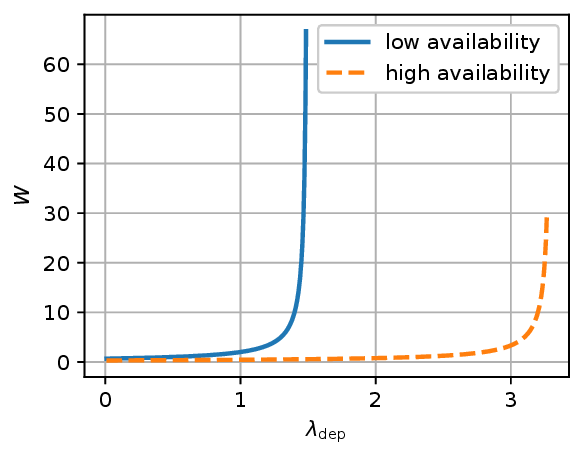}
        \caption{$W$ vs.\ $\lambda_\mathrm{dep}$}
        \label{fig:markov_operation_c}
    \end{subfigure}
    \caption{
        Operational implications of the four-state model.  
        (a) shows the effective service rate $\mu_{\mathrm{eff}}$, obtained from the stationary open-state probability $\pi_3$, as a function of $P_{\mathrm{good}}$.  
        (b) visualizes $\pi_3$ over the $(P_{\mathrm{good}}, P_{\mathrm{accept}})$ space, with a dashed contour indicating the stability boundary defined by $\lambda_\mathrm{dep} = c\,\pi_3$.  
        (c) illustrates the resulting delay behavior as demand $\lambda_\mathrm{dep}$ approaches system capacity, highlighting the sharp rise in delay near the stability threshold.  
    }
    \label{fig:markov_operation}
\end{figure*}
We use the steady-state probabilities of the Markov chain as a compact representation of operational performance, relating long-run state occupancy to effective service capacity, throughput, and delay. This interpretation clarifies how the parameters \(P_{\mathrm{good}}\), \(P_{\mathrm{accept}}\), and \(P_{\mathrm{success}}\) influence system behavior. Let \(c\) denote the maximum number of departures that can be served in each decision period when the fix is in the \emph{Gate Opened} (State~3), and let \(\lambda\) denote the average arrival rate of ready-to-depart flights. To translate the Markov chain steady-state behavior into operational performance measures, we approximate the departure process using standard queueing models.

\paragraph{Effective service rate}
Since departures can only be processed when the system is in State~3, and State~3 is occupied with long-run probability \(\pi_3\), the effective service rate is $\mu_{\mathrm{eff}} \triangleq c \,\pi_3.$ This quantity aggregates the effects of weather and decision-making into a single scalar that characterizes the usable departure capacity of the system.

\paragraph{Stability condition}
A stylized stability condition can be obtained by viewing the departure process as a single-server queue with arrival rate \(\lambda_\mathrm{dep}\) and service rate \(\mu_{\mathrm{eff}}\). In this abstraction, the system remains stable if the effective service rate exceeds the arrival rate, which is $\lambda_\mathrm{dep} \;<\; \mu_{\mathrm{eff}} \;=\; c\,\pi_3.$ This inequality defines a stability boundary in the \((P_{\mathrm{good}}, P_{\mathrm{accept}}, P_{\mathrm{success}})\) space through the dependence of \(\pi_3\) on the transition probabilities.

\paragraph{Approximate delay analysis}
Under a standard \(M/M/1\) approximation with arrival rate \(\lambda_\mathrm{dep}\) and service rate \(\mu_{\mathrm{eff}}\), the mean waiting time in the system satisfies
\begin{equation}
W \;\approx\; \frac{1}{\mu_{\mathrm{eff}} - \lambda_\mathrm{dep}}
\quad\text{for } \lambda_\mathrm{dep} < \mu_{\mathrm{eff}},
\end{equation}
and, by Little's law \cite{doi:10.1287/opre.9.3.383}, the corresponding mean queue length is
\begin{equation}
L \;\approx\; \lambda_\mathrm{dep}\, W \;=\; \frac{\lambda_\mathrm{dep}}{\mu_{\mathrm{eff}} - \lambda_\mathrm{dep}}.
\end{equation}
These expressions show how improvements in \(P_{\mathrm{good}}, P_{\mathrm{accept}},\) and \(P_{\mathrm{success}}\), which increase \(\pi_3\) and therefore \(\mu_{\mathrm{eff}}\), propagate directly to reductions in delay. In particular, when the system operates close to the stability boundary, small increases in \(\pi_3\) can yield disproportionately large reductions in expected waiting time and queue length.

The key relationships are summarized visually in \Cref{fig:markov_operation}, where all numerical examples use a service capacity of $c = 6$. The figure provides three complementary views of how the Markov chain is used to analyze operational performance. \Cref{fig:markov_operation_a} shows how the stationary availability $\pi_3$ determines the effective service rate $\mu_{\mathrm{eff}} = c \pi_3$. As $P_{\mathrm{good}}$ increases, $\pi_3$ rises more quickly, so the effective service rate improves at an accelerating pace. \Cref{fig:markov_operation_b} displays the stability boundary in the $(P_{\mathrm{good}}, P_{\mathrm{accept}})$ space where the condition $\lambda_\mathrm{dep} < c \pi_3$ holds. \Cref{fig:markov_operation_c} illustrates how delay responds to increasing demand using two representative availability levels, $\pi_3 = 0.25$ and $\pi_3 = 0.55$. These values show how higher availability substantially reduces delay near the stability boundary.

\subsection{Worst Case Scenario Evaluation}\label{sec:worst_result}
\subsubsection{Baseline Analysis}\label{sec:case1_result}
We begin with a baseline example that illustrates how the worst case analysis can be interpreted under independent flight behavior.
\begin{figure}[htbp]
    \centering
    \includegraphics[width=0.7\linewidth]{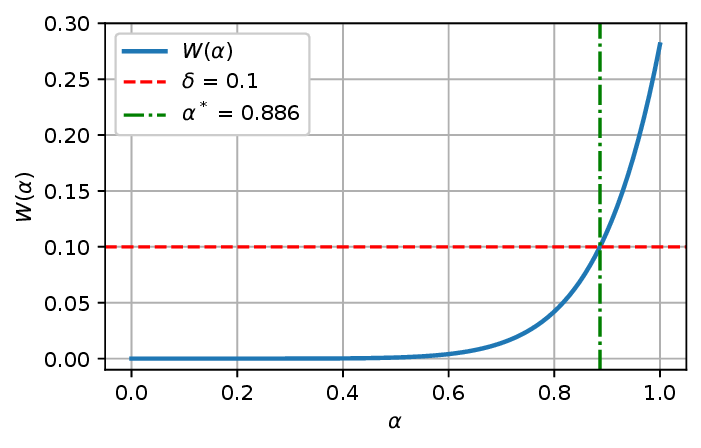}
    \caption{Worst case probability $W(\alpha)$ as a function of the rejective flight ratio $\alpha$. The red dashed horizontal line indicates the system failure threshold $\delta = 0.1$, and the green dash-dotted vertical line shows the critical point $\alpha^*$ where $W(\alpha^*) = \delta$, under $n = 10$, $U^- = -2$, $U^+ = 2$, and $\beta = 1$.}
    \label{fig:vanila_worst}
\end{figure}
The behavior of the worst case system failure probability $W(\alpha)$, defined as the probability of complete rejection when a fraction $\alpha$ of flights are rejective, and the associated tipping point $\alpha^*$ are illustrated in \Cref{fig:vanila_worst} under a representative set of parameters. 

This result provides a clear baseline for assessing system vulnerability when flights act independently. By identifying the tipping point beyond which widespread rejection triggers system failure, operators can proactively adjust pathfinder selection policies to remain below this risk boundary. This analysis shows that system vulnerability can be evaluated in a simple and interpretable way through the worst case failure probability, without requiring detailed assumptions on coordination or dynamics. For example, as illustrated in \Cref{fig:vanila_worst}, if the tolerable probability of complete rejection is set to $\delta = 0.1$, the corresponding threshold $\alpha^*$ implies that no more than 88.6\% of flights can be rejective before the system enters the fragile regime, i.e., when recovery is unlikely.

\subsubsection{Impact of Flight Selflessness}\label{sec:case2_result}
We next extend the baseline analysis by allowing flights to exhibit selfless behavior, whereby individual acceptance decisions respond to the observed level of collective rejection.
\begin{figure*}[htbp]
    \centering
    \begin{subfigure}{0.32\linewidth}
        \centering
        \includegraphics[width=\linewidth]{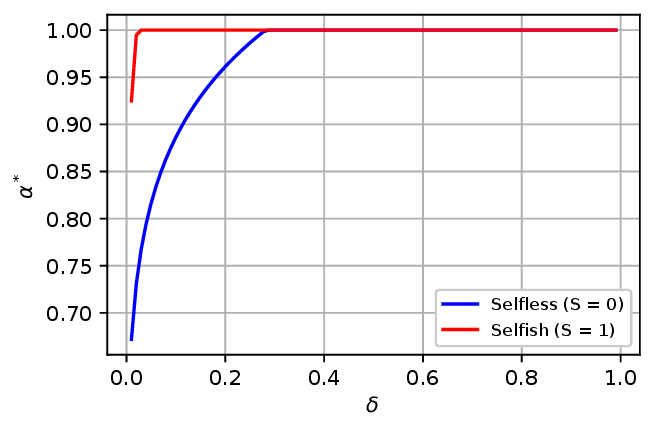}
        \caption{$\alpha^*$ vs.\ $\delta$}
        \label{fig:social_worst_a}
    \end{subfigure}
    \hfill
    \begin{subfigure}{0.32\linewidth}
        \centering
        \includegraphics[width=\linewidth]{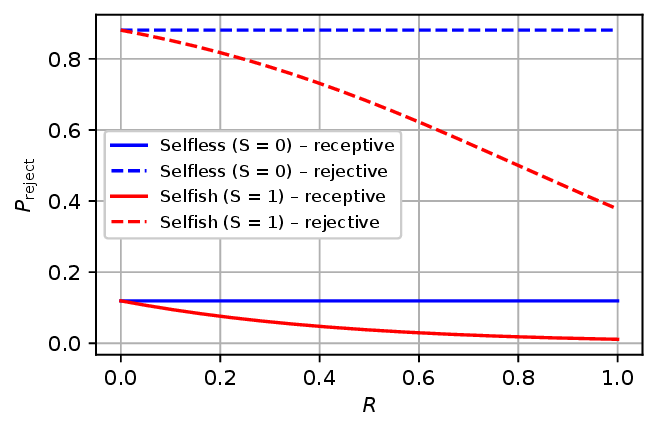}
        \caption{$P_{\mathrm{reject}}$ vs.\ $R$}
        \label{fig:social_worst_b}
    \end{subfigure}
    \hfill
    \begin{subfigure}{0.32\linewidth}
        \centering
        \includegraphics[width=\linewidth]{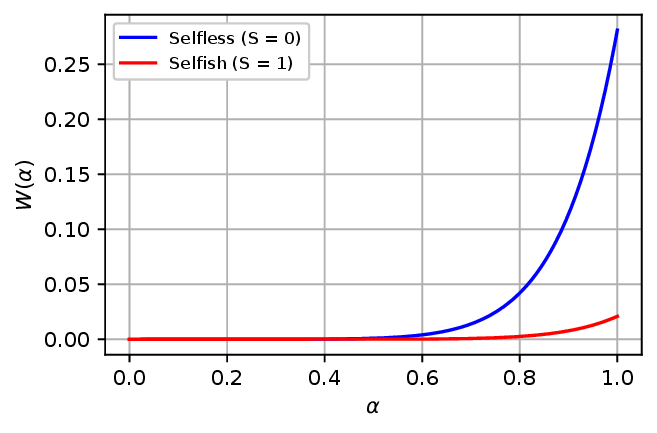}
        \caption{$W(\alpha)$ vs.\ $\alpha$}
        \label{fig:social_worst_c}
    \end{subfigure}

    \caption{
        Worst case collective failure analysis under varying degrees of selfishness 
        ($S=0$ for selfless, $S=1$ for selfish flights) with parameters 
        $n=10$, $U^+=2$, $U^-=-2$, $\beta=1$, $\gamma=2.5$, and $R=0.5$.  
        (a) illustrates how the tipping point $\alpha^*$ varies with the system failure tolerance $\delta$.  
        (b) reports the flight rejection probability $P_{\mathrm{reject}}$ as a function of the observed collective rejection ratio $R$, where solid and dashed lines indicate receptive and rejective groups.  
        (c) presents the worst case system failure probability $W(\alpha)$ across different values of $\alpha$.  
    }
    \label{fig:social_worst}
\end{figure*}
\Cref{fig:social_worst} illustrates how selfless behavior improves system robustness. As shown in \Cref{fig:social_worst_a}, the critical threshold \( \alpha^* \), at which the system failure probability \( W(\alpha) \) exceeds the tolerance level \( \delta \), increases as flights become more selfless. When flights are fully selfless (\( S=0 \)), the system can tolerate a much larger fraction of rejective flights before failure, and \( \alpha^* \) approaches one even for small values of \( \delta \), indicating strong system resilience. In contrast, when flights are selfish, \( \alpha^* \) grows slowly with \( \delta \) and reaches high values only under relatively loose failure tolerances. 

\Cref{fig:social_worst_b} further shows that rejection probabilities respond to the observed rejection rate \( R \) only under selfless behavior: as \( R \) increases, selfless flights become increasingly less likely to reject. Finally, \Cref{fig:social_worst_c} confirms that the overall worst case failure probability \( W(\alpha) \) is substantially lower under selfless behavior across all values of \( \alpha \). Together, these results demonstrate that even limited degrees of selflessness at the flight level can substantially enhance system resilience, providing operators with quantitative insight into how individual cooperation scales into overall operational stability.

\subsubsection{Impact of Environmental Uncertainty}\label{sec:case3_result}
We further extend the worst case analysis by introducing environmental uncertainty, which perturbs flights’ perceived utilities and acceptance decisions.
\begin{figure}[htbp]
    \centering
    \begin{subfigure}{0.49\linewidth}
        \centering
        \includegraphics[width=\linewidth]{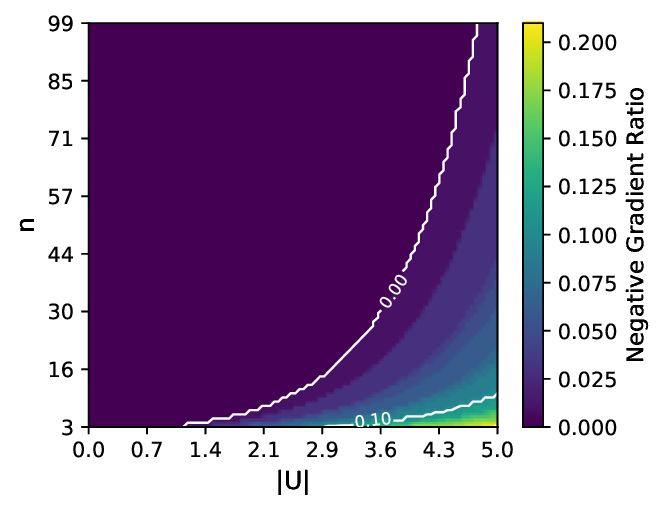}
        \caption{Rademacher}
        \label{fig:external_worst_rademacher}
    \end{subfigure}
    \hfill
    \begin{subfigure}{0.49\linewidth}
        \centering
        \includegraphics[width=\linewidth]{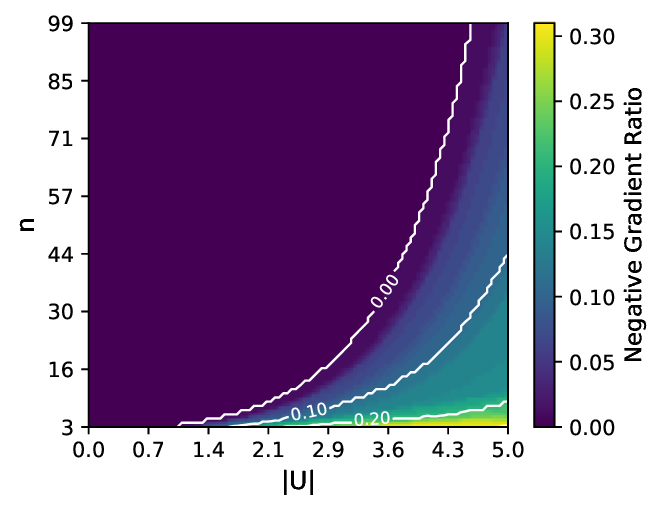}
        \caption{Gaussian}
        \label{fig:external_worst_gaussian}
    \end{subfigure}
    \caption{
    Comparison of gradient behavior under two uncertainty models.
    (a) shows the fraction of $(\alpha,\kappa)$ pairs with $\partial W / \partial \kappa < 0$ under Rademacher noise.
    (b) shows the corresponding fraction for Gaussian noise, with $\sigma$ denoting the uncertainty scale.
    Contour lines indicate equal gradient levels.
    }
    \label{fig:external_worst}
\end{figure}
\Cref{fig:external_worst} depicts the partial derivative $\partial W / \partial \theta$, which quantifies the effect of increasing environmental uncertainty on the worst case probability, as a function of both the rejective flight proportion $\alpha$ and the noise scale $\theta$. Here, $\theta$ represents the noise magnitude, with $\theta = \kappa$ for the Rademacher distribution and $\theta = \sigma$ for the Gaussian distribution. The gradient is evaluated over a grid of $(\alpha, \theta)$ values for each combination of flight count $n$ and absolute utility level $|U|$, with $\alpha \in [0,1]$ and $\theta \in [0,10]$.

The sign of $\partial W / \partial \theta$ indicates whether environmental randomness stabilizes or destabilizes the system. A positive value, $\partial W / \partial \theta > 0$, implies that increasing noise raises the likelihood of complete rejection and weakens system stability. In contrast, a negative value, $\partial W / \partial \theta < 0$, indicates that uncertainty suppresses the failure probability and enhances system robustness.

Across both model variants shown in \Cref{fig:external_worst}, negative values of $\partial W / \partial \theta$ are observed only rarely over most combinations of absolute utility $|U|$ and flight count $n$. Even in cases where negative gradients appear across the full $(\alpha,\theta)$ domain, the proportion of parameter settings for which $\partial W / \partial \theta < 0$ remains small. These results indicate that, in the majority of scenarios, increasing uncertainty $\theta$ tends to elevate the risk of complete rejection, thereby rendering the system more fragile.

Nevertheless, localized regions exist in which the gradient becomes negative, particularly when the absolute utility magnitude \( |U| \) is large and the number of flights \( n \) is small. Such conditions correspond to situations with unusually strong incentives or penalties (e.g., urgent departures or high operational costs) and limited pathfinder availability due to low traffic or constrained airspace. The Gaussian noise case exhibits slightly more frequent negative gradients than the Rademacher case. Across the simulation results, negative gradients are observed primarily when the rejective flight ratio \( \alpha \) is large, meaning most flights are already inclined to reject. 

In these cases, added randomness can cause a small number of flights to deviate from rejection, thereby preventing complete system failure. These results indicate that noise can play a stabilizing role under highly adverse rejection conditions. Overall, while uncertainty generally destabilizes the system in our model, these edge cases highlight conditions under which randomness can improve operational resilience within the considered setting, offering useful insight for risk-aware decision-making under uncertainty.

\subsection{Pathfinding Sequence Optimization Results}\label{sec:sequence_result}
\subsubsection{Data and Parameters Specification}
\label{sec:dataset}
\begin{table*}[htbp]
\centering
\caption{JFK departure dataset used for simulation. Pathfinder candidate flights are shown in bold.}
\label{tab:jfk_data}
\begin{tabular}{l l l l l r c}
\hline
\textbf{Flight Index} & \textbf{Aircraft} & \textbf{Destination} & 
\textbf{Scheduled Dep.} & \textbf{Scheduled Arr.} &
\textbf{Flight Time (min)} & \textbf{Wake Category} \\
\hline
\textbf{DAL1}   & \textbf{B764} & \textbf{London Heathrow (LHR)} & \textbf{21:12 EDT} & \textbf{08:41 BST} & \textbf{395} & \textbf{H} \\
JBU641 & A320 & Fort Lauderdale (KFLL) & 21:13 EDT & 23:53 EDT & 170 & M \\
RPA4535 & E75L & Indianapolis (KIND) & 21:14 EDT & 23:28 EDT & 109 & S \\
\textbf{SWR15}  & \textbf{A333} & \textbf{Zurich (ZRH)} & \textbf{21:15 EDT} & \textbf{10:22 CEST} & \textbf{433} & \textbf{H} \\
UPS2924 & B763 & Louisville (KSDF) & 21:18 EDT & 23:38 EDT & 140 & H \\
\textbf{DAL100} & \textbf{B764} & \textbf{Geneva (GVA)} & \textbf{21:20 EDT} & \textbf{10:21 CEST} & \textbf{427} & \textbf{H} \\
\textbf{THY2}   & \textbf{B77W} & \textbf{Istanbul (IST)} & \textbf{21:21 EDT} & \textbf{13:17 TRT} & \textbf{542} & \textbf{H} \\
\textbf{AAL292} & \textbf{B789} & \textbf{Delhi (DEL)} & \textbf{21:23 EDT} & \textbf{20:30 IST} & \textbf{823} & \textbf{H} \\
AAL1060 & A319 & Miami (KMIA) & 21:25 EDT & 00:28 EDT & 183 & M \\
JBU603  & A320 & Boston (KBOS) & 21:26 EDT & 22:13 EDT & 47 & M \\
\textbf{DAL52}  & \textbf{B763} & \textbf{Zurich (ZRH)} & \textbf{21:27 EDT} & \textbf{10:31 CEST} & \textbf{430} & \textbf{H} \\
\textbf{BAW172} & \textbf{B772} & \textbf{London Heathrow (LHR)} & \textbf{21:28 EDT} & \textbf{08:32 BST} & \textbf{370} & \textbf{H} \\
\textbf{SIA25} & \textbf{B77W} & \textbf{Frankfurt (FRA)} & \textbf{21:30 EDT} & \textbf{09:58 CEST} & \textbf{394} & \textbf{H} \\
UAL2019 & B738 & Charlotte (KCLT) & 21:53 EDT & 23:40 EDT & 107 & M \\
\textbf{JBU7}   & \textbf{A21N} & \textbf{London Heathrow (LHR)} & \textbf{21:56 EDT} & \textbf{09:10 BST} & \textbf{381} & \textbf{M} \\
\textbf{EIN106} & \textbf{A333} & \textbf{Dublin (DUB)} & \textbf{21:57 EDT} & \textbf{08:24 IST} & \textbf{64} & \textbf{H} \\
\textbf{ITY611} & \textbf{A339} & \textbf{Rome (FCO)} & \textbf{22:00 EDT} & \textbf{11:53 CEST} & \textbf{478} & \textbf{H} \\
JBU1991 & A320 & Rochester (KROC) & 22:01 EDT & 23:14 EDT & 73 & M \\
\textbf{RAM201} & \textbf{B788} & \textbf{Casablanca (CMN)} & \textbf{22:02 EDT} & \textbf{09:33 WEST} & \textbf{395} & \textbf{H} \\
DAL2920 & B739 & Indianapolis (KIND) & 22:03 EDT & 00:17 EDT & 134 & M \\
JBU1722 & A320 & Palm Beach (KPBI) & 22:04 EDT & 00:58 EDT & 174 & M \\
\textbf{JBU73} & \textbf{A21N} & \textbf{Edinburgh (EDI)} & \textbf{22:04 EDT} & \textbf{09:34 BST} & \textbf{393} & \textbf{M} \\
\textbf{IBE326} & \textbf{A333} & \textbf{Madrid (MAD)} & \textbf{22:06 EDT} & \textbf{11:08 CEST} & \textbf{425} & \textbf{H} \\
RPA4573 & E75L & Raleigh--Durham (KRDU) & 22:07 EDT & 23:22 EDT & 75 & S \\
RPA4659 & E75L & Boston (KBOS) & 22:08 EDT & 22:30 EDT & 22 & S \\
JBU102  & BCS3 & Buffalo (KBUF) & 22:09 EDT & 23:19 EDT & 70 & M \\
ASA41   & B739 & San Francisco (KSFO) & 22:11 EDT & 00:33 PDT & 322 & M \\
RPA4546 & E75L & Reagan National (KDCA) & 22:11 EDT & 23:03 EDT & 52 & S \\
\hline
\end{tabular}
\end{table*}
We use a historical set of departure schedules from New York John F. Kennedy International Airport (JFK) on September 23, 2025. The dataset, derived from FlightAware \cite{flightaware}, contains 30 departing flights scheduled between 22:00--23:30 EDT. All timestamps were converted to UTC and normalized by subtracting the earliest departure time. Aircraft wake categories and associated separation standards were assigned with reference to FAA wake turbulence procedures and recategorization guidance \cite{FAA711065}. Each aircraft in the dataset is classified into Small (S), Medium (M), or Heavy (H) categories based on certified maximum takeoff weight. Scheduled block times were calculated by taking the difference between the scheduled departure time (in EDT for JFK) and the scheduled arrival time expressed in the destination airport’s local time zone.

To construct the pathfinder candidate set, we filtered flights destined for eastbound locations, reflecting the operational convention that pathfinders are selected from aircraft routed through the same directional fix being assessed. This filtering resulted in a set of 14 candidate flights. In the absence of airline connection data, we assume the pathfinder's departure-time reduction benefit \(T_i\) equals its airline benefit \(B_{\mathrm{dep},i}\), a standard simplification when connection penalties are unavailable. \Cref{tab:jfk_data} summarizes all flights used in the simulation. Eastbound pathfinder candidate flights are highlighted in bold.

\subsubsection{Parameter Estimation via Discrete Event Simulation}\label{sec:des_result}
To obtain the numerical inputs required for the sequencing formulations, we apply the discrete-event simulator to the calibrated JFK departure environment. The simulation uses Runway~4L and Runway~31L, one of the primary parallel departure runway configurations at JFK \cite{PANYNJ_JFKAC_2022}. Wake-based separations and fix availability follow the operational settings introduced earlier. All flights are simulated under identical stochastic taxi-out and cancellation settings, with taxi-out times drawn from the same empirical distribution, so that differences in outcomes arise solely from the pathfinder intervention.

For each candidate flight \(i\) and each offer position \(k\), two simulations are executed with identical random seeds: (i) a \emph{baseline} run in which no flight receives pathfinder privileges, and (ii) a \emph{pathfinder} run in which flight \(i\) accepts the offer at slot \(k\). In the pathfinder run, offer-related overheads are applied, the designated flight is placed in the runway queue to achieve the earliest predicted departure, and the relevant fix is opened at the corresponding trigger time. Since both runs share the same underlying randomness, the paired outputs isolate the operational impact of granting pathfinder status to flight \(i\). The baseline–pathfinder comparisons yield, for each \((i,k)\), the five evaluation metrics \(T, B^{\text{dep}}, D^{\text{sys}}, G^\mathrm{ATC}, G^\mathrm{disp}\). 

\begin{figure*}[t]
    \centering
    \begin{subfigure}{0.195\textwidth}
        \centering
        \includegraphics[width=\linewidth]{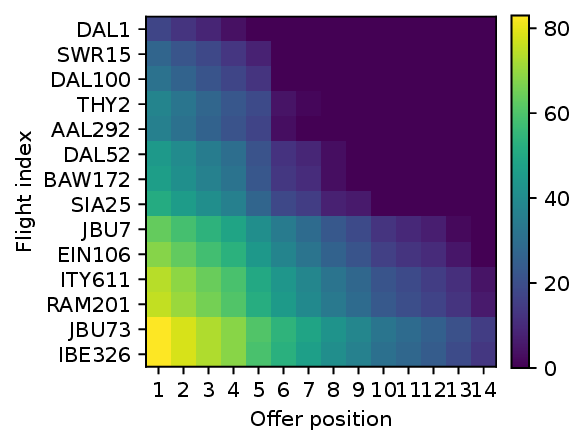}
        \caption{$T_{i,k}$}
    \end{subfigure}
    \hfill
    \begin{subfigure}{0.195\textwidth}
        \centering
        \includegraphics[width=\linewidth]{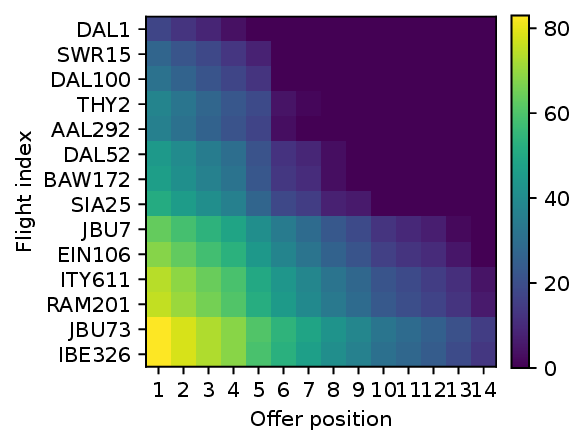}
        \caption{$B^{\mathrm{dep}}_{i,k}$}
    \end{subfigure}
    \hfill
    \begin{subfigure}{0.195\textwidth}
        \centering
        \includegraphics[width=\linewidth]{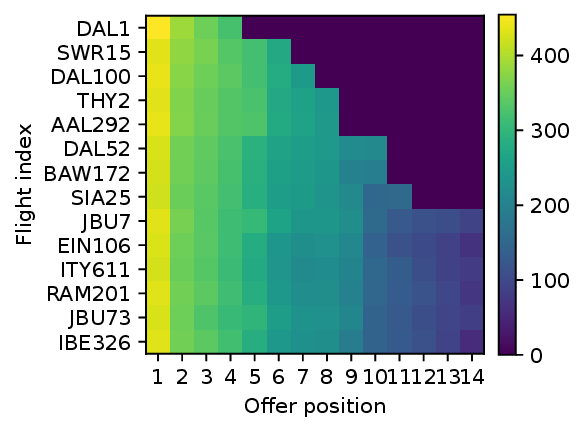}
        \caption{$D^{\mathrm{sys}}_{i,k}$}
    \end{subfigure}
    \hfill
    \begin{subfigure}{0.195\textwidth}
        \centering
        \includegraphics[width=\linewidth]{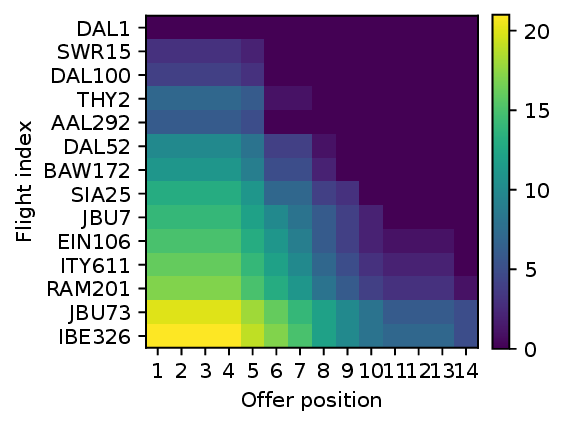}
        \caption{$G^{\mathrm{disp}}_{i,k}$}
    \end{subfigure}
    \hfill
    \begin{subfigure}{0.195\textwidth}
        \centering
        \includegraphics[width=\linewidth]{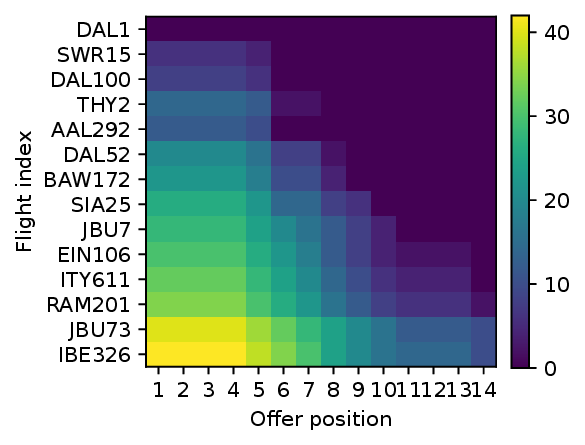}
        \caption{$G^{\mathrm{ATC}}_{i,k}$}
    \end{subfigure}
    \caption{
        Heatmaps of the five simulation-derived parameter matrices. Each matrix is indexed by candidate flight \(i\) and offer position \(k\), forming the sensitivity inputs for downstream sequencing optimization.
    }
    \label{fig:des_result}
\end{figure*}

The resulting parameter matrices are illustrated in \Cref{fig:des_result}. Candidate flights are indexed along the vertical axis according to their original scheduled departure order, with earlier scheduled flights shown toward the top and later flights toward the bottom. The parameters $T_{i,k}$, $B^{\mathrm{dep}}_{i,k}$, $G^{\mathrm{ATC}}_{i,k}$, and $G^{\mathrm{disp}}_{i,k}$ exhibit similar qualitative trends. Specifically, their values increase as the offer position moves earlier in the sequence and as the candidate flight’s original scheduled position moves later in the departure queue. This indicates that later-scheduled flights benefit more from early pathfinder designation due to larger queue-jumping opportunities and greater downstream delay reductions. The airline benefit $B^{\mathrm{dep}}_{i,k}$ would be further amplified if detailed passenger connection and aircraft rotation information were incorporated.

In contrast, the system delay reduction $D^{\mathrm{sys}}_{i,k}$ shows a much stronger dependence on the offer position than on the flight index. Larger values of $D^{\mathrm{sys}}_{i,k}$ consistently occur when the offer position is early, since early acceptance enables earlier fix opening and greater cumulative delay savings for downstream flights. For flights that are already near the front of the original departure schedule, however, $D^{\mathrm{sys}}_{i,k}$ drops to zero at later offer positions, as accepting the pathfinder attempt would place the flight behind its original queue position. This eliminates both the incentive to accept and the associated system benefit.

\subsubsection{Optimization Outcome Analysis}\label{sec:opt_result}
The ATC-initiated and dispatcher-initiated optimization formulations were evaluated using the calculated parameter matrices described in the previous subsection. For both formulations, all input matrices were normalized to the range $[0,1]$, and the acceptance probability was computed using a logistic model with a fixed success probability $P_{\mathrm{succ}}=0.9$ and a common sensitivity parameter $\beta$. The participation and failure costs were set to zero ($c_i=d_i=0$) in order to isolate the effect of departure-time reward on acceptance behavior, as reliable operational data for these cost components were not available. All combinations of the offer budget $B$, the risk penalty weight $\lambda$ (with $\lambda_{\mathrm{ATC}}$ and $\lambda_{\mathrm{disp}}$ denoting the ATC- and dispatcher-initiated instances), and the sensitivity $\beta$ over their respective ranges were evaluated. We used $B \in [3,12]$ (10 grid points), $\lambda \in [0,1]$ (11 grid points), and $\beta \in [0,5]$ (6 grid points), yielding a total of $10 \times 11 \times 6 = 660$ instances. For each formulation, a total of $660$ mixed-integer nonlinear programming instances were solved using the Bonmin solver.

To enable a direct structural comparison between ATC-driven and dispatcher-driven sequencing policies, both optimization formulations were evaluated over the same set of 14 candidate flights. This controlled setup allows the resulting sequencing behaviors to be attributed to the objective structure rather than to differences in candidate availability. In practice, dispatcher-initiated optimization would generally be restricted to flights belonging to the corresponding airline. 

\begin{table}[htbp]
\centering
\caption{Runtime summary of the sequence optimization.}
\label{tab:runtime_summary}
\begin{tabular}{c|cccc}
\hline
\textbf{Model} & \textbf{Avg. (s)} & \textbf{Std. (s)} & \textbf{Min (s)} & \textbf{Max (s)} \\
\hline
ATC & 14.26 & 28.92 & 0.28 & 233.41 \\
Dispatcher & 20.42 & 26.35 & 0.45 & 194.06 \\
\hline
\end{tabular}
\end{table}

The computational performance of the two optimization formulations is summarized in \Cref{tab:runtime_summary}. All instances were solved on a laptop equipped with a 13th Gen Intel(R) Core i7-1360P processor. The ATC formulation achieved an average runtime of $14.26$~s, while the dispatcher formulation required $20.42$~s on average. Even the longest solution time remained below four minutes for both models, demonstrating the practical feasibility of deploying the proposed pathfinder offer sequence optimization in near real-time operational settings.

\begin{figure}[htbp]
    \centering
    \begin{subfigure}{0.492\linewidth}
        \centering
        \includegraphics[width=\linewidth]{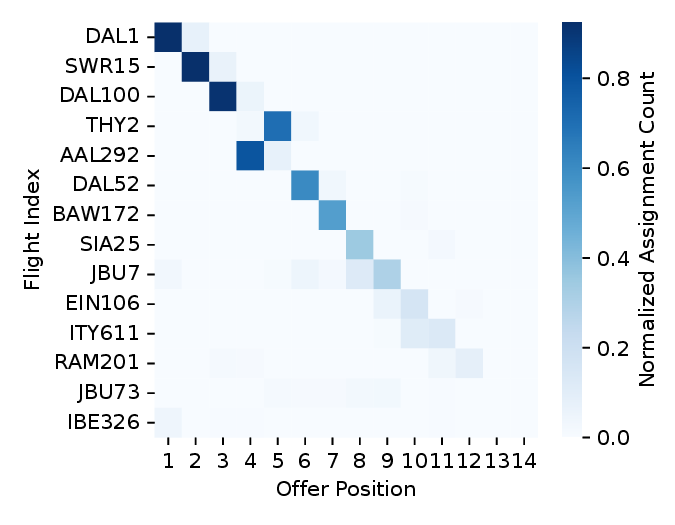}
        \caption{ATC}
        \label{fig:optimization_atc}
    \end{subfigure}
    \hfill
    \begin{subfigure}{0.492\linewidth}
        \centering
        \includegraphics[width=\linewidth]{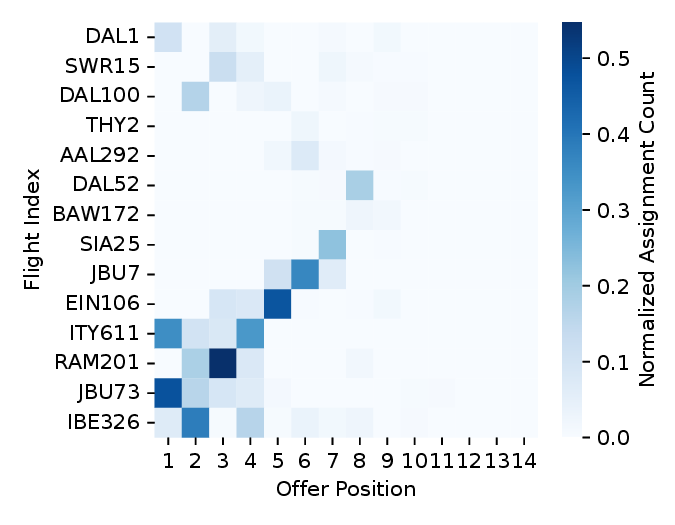}
        \caption{Dispatcher}
        \label{fig:optimization_disp}
    \end{subfigure}
    \caption{
    Comparison of optimal flight-to-offer assignment patterns under the ATC-initiated and dispatcher-initiated sequencing policies.
    Each cell indicates whether flight $i$ is assigned to offer position $k$.
    }
    \label{fig:opt_result}
\end{figure}
\paragraph{Contrasting sequencing structures under ATC and dispatcher objectives}
\Cref{fig:opt_result} visualizes representative optimal flight-to-offer assignment patterns obtained under the ATC-initiated and dispatcher-initiated sequencing policies. Clear structural differences emerge between the two sequencing policies. As shown in \Cref{fig:opt_result}a, the ATC formulation, which maximizes system delay reduction, tends to assign earlier-scheduled flights to the first few offer positions. In our simulation, $D^{\mathrm{sys}}_{i,k}$ depends primarily on the offer position $k$, while the ATC-side penalty $G^{\mathrm{ATC}}_{i,k}$ increases when very late flights are moved far forward in the queue. As a result, the ATC-driven sequences favor flights that are already near the front of the departure queue, capturing most of the available system benefit from early offers while avoiding large operational disruptions.

In contrast, the dispatcher formulation in \Cref{fig:opt_result}b more aggressively promotes later-scheduled flights to early offer positions. From the dispatcher’s perspective, the airline-specific benefit $B^{\mathrm{dep}}_{i,k}$ grows with how much a flight’s position in the departure queue is advanced, so moving flights from the back of the queue forward yields substantial local gains. Even when this behavior yields smaller system-wide delay reduction, the dispatcher objective still prefers such assignments as long as the perceived dispatcher-side risk $G^{\mathrm{disp}}_{i,k}$ remains acceptable. Since both policies are evaluated over the same candidate set, these contrasting patterns can be attributed to differences in the underlying objective structures rather than to differences in candidate availability.

\paragraph{Concentration of benefit in early offer positions}
Across the explored parameter settings, the expected objective value is dominated by the first few offer positions. Although the budget $B$ limits the maximum number of offers that can be issued, the contribution of later offers diminishes rapidly due to the stopping probability structure in the objective, which multiplies each offer’s benefit by the probability that all earlier offers are rejected. This concentration effect is a direct consequence of the multiplicative survival term, which attenuates the expected contribution of later offer positions. In our numerical experiments, the majority of the expected benefit is concentrated in the early part of the sequence, highlighting the importance of early acceptance in driving overall system performance for the considered parameter ranges.
\begin{figure}[htbp]
    \centering
    \begin{subfigure}{0.49\linewidth}
        \centering
        \includegraphics[width=\linewidth]{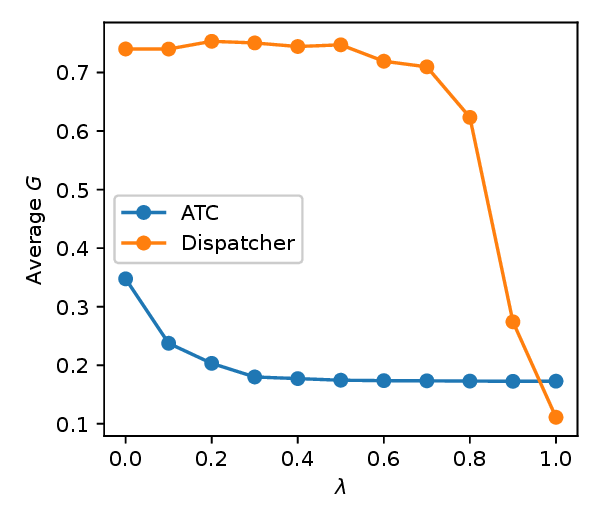}
        \caption{$\lambda$ sensitivity}
        \label{fig:lambda_sensitivity}
    \end{subfigure}
    \hfill
    \begin{subfigure}{0.49\linewidth}
        \centering
        \includegraphics[width=\linewidth]{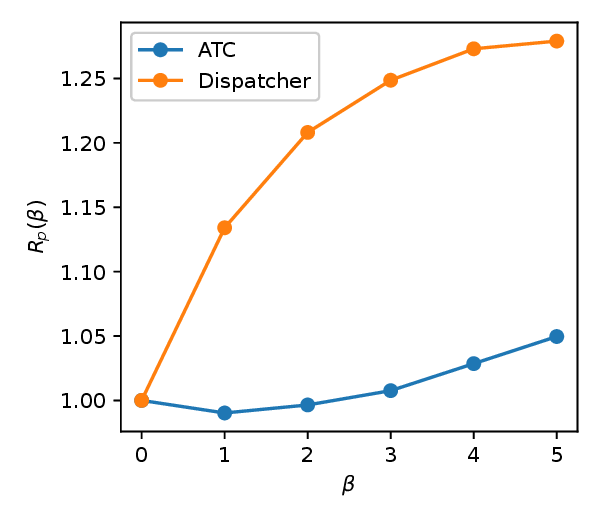}
        \caption{$\beta$ sensitivity}
        \label{fig:beta_sensitivity}
    \end{subfigure}
    \caption{Sensitivity analysis with respect to the risk penalty weight $\lambda$ and the acceptance sensitivity parameter $\beta$. 
    (a) shows the variation of the average risk metric $G$ of the selected flights as $\lambda$ increases.
    (b) depicts the relative selection ratio $R_p(\beta)$, defined as the ratio between the average acceptance probability of selected offers and that of the overall candidate set.}
    \label{fig:sensitivity_analysis}
\end{figure}
\paragraph{Risk sensitivity under ATC and dispatcher control}
\Cref{fig:lambda_sensitivity} examines how the sequencing policies respond to the risk‐penalty weight $\lambda$, which governs when operational risk begins to outweigh the benefit of assigning a flight as the pathfinder. Rather than inducing gradual adjustments, varying $\lambda$ produces threshold‐like shifts in sequencing behavior. Under the ATC formulation, most structural change occurs at low values of $\lambda \in [0,0.3]$, indicating that ATC sequencing becomes more risk-averse even under modest increases in risk, because system delay benefits diminish rapidly as flights are moved far forward in the queue. This behavior is consistent with interviews indicating that ATC decision-making prioritizes safety and conservative risk management \cite{interview}. In contrast, under the dispatcher formulation, sequencing behavior remains largely unchanged until $\lambda$ approaches the upper end of its range $\lambda \in [0.7,1]$, at which point the structure shifts sharply, reflecting the strong local incentives to advance late flights. These contrasting patterns reflect fundamentally different thresholds at which risk considerations begin to dominate sequencing decisions.

\paragraph{Sensitivity to acceptance likelihood}
To characterize the structural impact of $\beta$, we define the relative selection ratio
\begin{equation}
R_p(\beta) = \frac{\bar{p}_{\mathrm{sel}}(\beta)}{\bar{p}_{\mathrm{all}}(\beta)},
\end{equation}
where $\bar{p}_{\mathrm{sel}}(\beta)$ is the average acceptance probability of the selected offers and $\bar{p}_{\mathrm{all}}(\beta)$ is the corresponding average over the entire candidate set. Larger values of $R_p$ indicate stronger preference for high-acceptance candidates. 

As shown in \Cref{fig:beta_sensitivity}, increasing $\beta$ leads both formulations to place greater emphasis on acceptance likelihood, resulting in higher values of $R_p$. This reflects a shift toward more reliability-driven sequencing, in which the optimizer prioritizes securing an early acceptance over potential delay or departure benefits. The dispatcher formulation exhibits a steeper increase in $R_p$, indicating that dispatcher-driven sequencing more aggressively filters for high-acceptance candidates. This behavior arises because the acceptance probability $p_{i,k}$ is reinforced by the aligned trends of the dispatcher’s benefit $B^{\mathrm{dep}}$ and perceived risk $G^{\mathrm{disp}}$, causing acceptance likelihood to play a more dominant role in the sequencing decision.

\paragraph{Implications for pathfinder sequencing design}
Together, the $\lambda$ and $\beta$ sensitivity analyses show that sequencing behavior is driven by both the magnitude of the policy parameters and by whose objective is used to order offers. Even under identical values of $\lambda$ and $\beta$, ATC- and dispatcher-initiated sequencing produce qualitatively different thresholds and selection patterns, leading to structurally different recovery dynamics. The choice of sequencing mechanism is therefore a first-order design decision for pathfinder operations, with implications for both system reliability and airline profits.

\subsection{Key Takeaways}\label{sec:key_takeaways}
We summarize the main insights from the numerical results.
\begin{itemize}
\item \textit{Acceptance probability drives performance:} 
Fix availability and delay vary significantly with $P_{\mathrm{accept}}$.

\item \textit{Fix availability governs capacity:} 
Steady-state fix availability determines throughput and stability.

\item \textit{Collective rejection is a primary failure mode:} 
Recovery can fail when many flights are unlikely to accept.

\item \textit{Modest selfless behavior improves resilience:} 
Limited selflessness increases tolerance to rejective populations.

\item \textit{Early offers concentrate recovery benefit:} 
Most benefit is realized in early offer positions.

\item \textit{Sequencing mechanism influences recovery dynamics:} 
ATC and dispatcher sequencing induce different offer patterns.
\end{itemize}

\section{Conclusion}\label{sec:conclusion}
This paper modeled pathfinder operations as a stochastic, multi-agent decision-making problem under convective weather. Using a Markov chain, we characterized system transitions among fix closure, pathfinder selection, execution, and reopening, and linked these dynamics to effective departure capacity, queue stability, and delay. We introduced utility-based flight acceptance models and showed how individual decisions can escalate into system vulnerability through collective rejection, while selfless behavior and environmental uncertainty were found to significantly shape system robustness. Building on these behavioral foundations, we developed optimization-based sequencing models for both ATC- and dispatcher-initiated pathfinder selection. Simulation experiments using real departure schedules and data from JFK demonstrated that the two institutional perspectives yield systematically different offer strategies, and that early acceptance dominates overall system benefit. Overall, this work provides a unified and practical framework linking weather uncertainty, individual incentives, sequential decision-making, and system performance, laying a foundation for future data-driven decision support for pathfinder operations.

\subsection{Limitations of Work and Future Directions}\label{sec:future_work}
Although this study presents a first formal analytical framework for pathfinder operations, several simplifying assumptions were made to preserve tractability, and naturally provides directions for future work. The Markov chain can be extended to nonstationary or regime-switching formulations that incorporate real-time weather forecasts, and the utility model can incorporate richer operational factors such as passenger connections, aircraft rotations, and crew constraints as data become available. The worst case rejection analysis can likewise be expanded to capture learning, coordination, and information sharing among airlines. Future work may also relax the homogeneous acceptance assumption by allowing different flights and airlines to have different behavioral sensitivities estimated from data. From an implementation perspective, scalable algorithms and additional operational considerations will be needed to support deployment in larger, multi-fix environments. While the simulation was calibrated to JFK, the framework itself is general and designed to support the development of practical, data-driven decision support tools for weather-impacted air traffic operations.

\section*{Acknowledgments} \label{sec:ack}

\noindent
We thank Curt Rademaker, Bill Bateman, Ron Foley, and Wayne Hubbard for helpful discussions and feedback on this work.

\section*{Disclaimer} \label{sec:disc}

\noindent
The contents of this document reflect the views of the authors and do not necessarily reflect the views of the National Aeronautics and Space Administration (NASA), the Federal Aviation Administration (FAA), or the Department of Transportation (DOT). Neither NASA, FAA, nor DOT make any warranty or guarantee, expressed or implied, concerning the content or accuracy of these views.

\bibliographystyle{IEEEtran}
\bibliography{reference}  

\end{document}